\documentclass[12pt]{amsart}

\usepackage{hyperref}
\usepackage{amssymb,amsmath}
\usepackage[alphabetic,nobysame,abbrev]{amsrefs}
\usepackage{enumerate}
\usepackage{color}
\usepackage{pdfsync}
\RequirePackage{fix-cm}
\usepackage{array}
\usepackage{longtable}
\usepackage{float}

\DeclareFontFamily{U}{mathx}{}
\DeclareFontShape{U}{mathx}{m}{n}{<-> mathx10}{}
\DeclareSymbolFont{mathx}{U}{mathx}{m}{n}
\newlength\savedwidth 

\usepackage{fancybox,framed}
\usepackage{bm}
\usepackage{mdframed}
\usepackage{cancel}
\usepackage{makecell}
\usepackage{pdfpages}
\usepackage{pgf,tikz}
\usepackage{mathrsfs}
\usetikzlibrary{arrows}
\usetikzlibrary{arrows,decorations.markings}
\usetikzlibrary{matrix}
\usetikzlibrary{backgrounds}

	\usepackage{diagbox}

\definecolor{qqqqff}{rgb}{0.,0.,1.}
\definecolor{xdxdff}{rgb}{0.49019607843137253,0.49019607843137253,1.}

\definecolor{qqqqff}{rgb}{0.,0.,1.}

\usepackage{tikz}
\usepackage{graphicx}
\usepackage[cmtip,all]{xy}

\usepackage[draft]{todonotes}
\usepackage[cmtip,all]{xy}

\usepackage{hyperref}
\usepackage[utf8]{inputenc}
\usepackage{amsmath,amsfonts,amssymb,euscript,graphicx,color,tikz}

\theoremstyle{plain}
\newtheorem{theorem}{Theorem}[subsection]
\newtheorem{thm}[theorem]{Theorem}
\newtheorem{lem}[theorem]{Lemma}
\newtheorem{cor}[theorem]{Corollary}
\newtheorem{pro}[theorem]{Proposition}

\theoremstyle{definition}
\newtheorem{DEF}[theorem]{Definition}

\newtheorem{exa}[theorem]{Example}
\newtheorem{con}[theorem]{Convention}
\newtheorem{const}[theorem]{Construction}

\newtheorem{rem}[theorem]{Remark}

\newtheorem{parag}[theorem]{{}}

\newcommand\pref[1]{\textbf{\ref{#1}}}
\numberwithin{equation}{section}

\newcommand{\La}{\langle}

\newcommand{\sub}{\subseteq}

\newcommand{\la}{Lie algebra }

\newcommand{\fm}{(\cdot,\cdot)}
\newcommand{\fh}{\mathfrak{h}}

\newcommand{\fg}{\mathfrak{g}}
\newcommand{\tfg}{\tilde{\mathfrak{g}}}

\newcommand{\ep}{\hfill$\Box$}

\def\ad{\hbox{ad}}
\def\andd{\quad\hbox{and}\quad}
\def\sg{\sigma}
\def\a{\alpha}
\def\b{\beta}
\def\lam{\lambda}
\def\Lam{\Lambda}
\def\ep{\epsilon}
\def\andd{\quad\hbox{and}\quad}

\def\id{\hbox{id}}
\def\Aut{\hbox{Aut}}

\def\g{\mathfrak{g}}

\def\mod{\hbox{mod}}
\def\andd{\quad\hbox{and}\quad}
\def\ind{\hbox{ind}}

\def\v{{\mathcal V}}
\def\u{{\mathcal U}}
\def\vd{\dot{\mathcal V}}

\def\fm{(\cdot,\cdot)}
\def\a{\alpha}

\def\sub{\subseteq}
\def\rd{\dot{R}}

\def\lam{\lambda}
\def\Lam{\Lambda}

\def\1k{\frac{1}{k}}

\def\la{\langle}
\def\ra{\rangle}
\def\rds{\dot{R}_{sh}}
\def\rdl{\dot{R}_{lg}}

\def\d{\delta}

\def\b{\beta}

\def\qed{\hfill$\Box$}

\def\sg{\sigma}

\def\hh{{\mathcal H}}

\def\sg{\sigma}

\usepackage{verbatim} 

\def\ad{\hbox{ad}}

\def\bbbz{{\mathbb Z}}

\def\bbbr{{\mathbb R}}

\def\bbbk{{\mathbb K}}

\def\zn{{\mathbb Z}^{\nu}}
\def\aa{\mathcal A}

\def\tr{\hbox{tr}}

\def\dd{\mathcal D}
\def\ep{\epsilon}

\def\aff{\hbox{Aff}}

\def\jj{{\mathcal J}}

\def\gl{\frak{g}\ell}

\def\proof{{\noindent\bf Proof. }}

\def\rds{\dot{R}_{sh}}
\def\rdl{\dot{R}_{lg}}

\def\g{\mathfrak{g}}

\def\span{\hbox{span}}

\def\cent{\hbox{Cent}}

\def\Hom{\hbox{Hom}}
\def\da{\dot\alpha}

\def\rank{\hbox{rank}}

\def\ss{\mathfrak{S}}

\def\DynkinNodeSize{1.5mm}
\def\DynkinArrowLength{2mm}
\tikzset{
	dnode/.style={
		circle,
		inner sep=0pt,
		minimum size=\DynkinNodeSize,
		fill=white,
		draw},
	middlearrow/.style={
		decoration={markings,
			mark=at position 0.8 with
			{\rdaw (0:0mm) -- +(+140:\DynkinArrowLength); \rdaw (0:0mm) -- +(-140:\DynkinArrowLength);},
		},
		postaction={decorate}
	},
	leftrightarrow/.style={
		decoration={markings,
			mark=at position 0.999 with
			{
				\rdaw (0:0mm) -- +(+135:\DynkinArrowLength); \rdaw (0:0mm) -- +(-135:\DynkinArrowLength);
			},
			mark=at position 0.001 with
			{
				\rdaw (0:0mm) -- +(+45:\DynkinArrowLength); \rdaw (0:0mm) -- +(-45:\DynkinArrowLength);
			},
		},
		postaction={decorate}
	},
	sedge/.style={
	},
	dedge/.style={
		middlearrow,
		double distance=0.6mm,
	},
	tedge/.style={
		middlearrow,
		double distance=1.0mm+\pgflinewidth,
		postaction={draw}, 
	},
	infedge/.style={
		leftrightarrow,
		double distance=0.5mm,
	},
}

\begin{document}

%
%
\title{Localization and filtration in  extended affine Lie algebras}

\author{S. Azam}
\address
{Department of Pure Mathematics\\Faculty of Mathematics and Statistics\\
	University of Isfahan\\ P.O.Box: 81746-73441\\ Isfahan, Iran, and\\
	School of Mathematics, Institute for
	Research in Fundamental Sciences (IPM), P.O. Box: 19395-5746.} \email{azam@ipm.ir, azam@sci.ui.ac.ir}
\thanks{This research was in part carried out in
	IPM-Isfahan Branch.}


\begin{abstract}
We investigate the notions of \emph{localization} and \emph{filtration} in the context of extended affine Lie algebras. Our primary objective is to develop a localization theory that facilitates the construction of meaningful local substructures, particularly local affine Lie subalgebras. These subalgebras play a crucial role in understanding the global structure of extended affine Lie algebras.  It is noteworthy that the existence of appropriate local subalgebras, particularly affine Lie subalgebras, is also fundamental to the representation theory of extended affine Lie algebras. As a natural outcome of our localization approach, we also introduce a formal notion of filtration for a given extended affine Lie algebra. This study is motivated by our interest in modular theory, specifically the integral structures of extended affine Lie algebras.
\end{abstract}


\maketitle

\section{Introduction}\label{intro}

Extended affine Lie algebras (EALAs) and extended affine root systems (EARSs) are higher-nullity generalizations of affine Lie algebras and root systems that preserve and extend many of their foundational structural properties. The class of EALAs was first introduced in 1990 under the name "quasi-simple Lie algebras" \cite{H-KT90}, and its theory was subsequently developed in 1997 under the now-standard terminology "extended affine Lie algebras" \cite{AABGP97}. The notion of EARSs was formally introduced earlier, in 1985, by K. Saito \cite{Sai85}.

 A distinctive characteristic of EALAs and EARSs is their axiomatic definitions, which permit a unified and rigorous treatment. As expected, the root system of an EALA is itself an EARS. In this work, we focus our attention on EALAs of reduced type.

Like their affine counterparts, an EALA possesses a non-degenerate symmetric invariant bilinear form, a root space decomposition relative to a self-centralizing Cartan subalgebra, and the root system satisfying specific finiteness and discreteness conditions.

Since their introduction, two principal approaches have shaped the study of EARSs: the embedding of local finite root subsystems and the embedding of local affine subsystems. For instance, Saito and collaborators \cite{Sai85} adopted the latter perspective, using local affine subsystems as the foundation for their classification of EARSs. In contrast, the authors of \cite{AABGP97} and subsequent works employed local finite subsystems. The finite perspective also appears in the study of associated Weyl groups, as one can see for example in \cite{MS92}.

A powerful method in Lie theory involves analyzing local substructures with well-understood properties. The study of $\mathfrak{sl}_2$-subalgebras has been particularly instrumental in uncovering the structure and representation theory of both finite and infinite-dimensional Lie algebras. This approach extends naturally to generalizations such as Kac-Moody and extended affine Lie algebras, as for example is illustrated in \cite[Chapter I]{AABGP97}.

Within extended affine Lie theory, both local and canonical subalgebras serve as essential analytical tools. Among canonical subalgebras, the \emph{core}, defined as the subalgebra generated by non-isotropic root spaces, constitutes the structural backbone of the entire EALA and has been central to the structure theory and classification of EALAs since their introduction.

This article develops a localization theory for EALAs, inspired by our previous works on modular theory and integral structures in EALAs. The natural starting point is to consider  a {\it closed} local subsystem of the root system of a given EALA, and look at the corresponding subalgebra. However, this subalgebra may fail to be an EALA. To address this, we introduce a method of cutting and pasting subspaces to obtain a (pure) extended affine subalgebra. Our construction of local extended affine Lie subalgebras reveals profound interactions with the core, establishing new connections between local and canonical substructures.

The central contribution of this work is the construction of local affine Lie subalgebras within a given extended affine Lie algebra of non-zero nullity, which again justifies the name EALAs for this class of Lie algebras. This construction provides a powerful tool for both structural analysis and representation theory, enabling classical techniques to be applied in localized settings. Our approach facilitates problem reduction to these well-understood subalgebras, followed by globalization of results. This strategy is inspired by and employed in our series of papers \cites{AF19,AFI22,AI23,Az24,Az25} on integral structures and Chevalley bases for EALAs.

Complementing our localization approach, we introduce a filtration technique that systematically decomposes the algebra based on a filtration of the underlying extended affine root system. As is well known from classical Lie theory, an increasing filtration makes a Lie algebra more tractable by replacing it with a graded version. In this spirit, we anticipate that the combined framework of localization and filtration will prove fruitful for analyzing the structure and representations of EALAs.

The paper is organized as follows. Section \ref{preliminaries} reviews essential definitions and foundational properties of extended affine root systems and Lie algebras. Section \ref{sec:tame6} briefly introduces E. Neher’s construction of EALAs \cite{Neh04}, which is particularly relevant for the examples presented in Section \ref{examples}.

In Section \ref{sec:local}, we consider an EALA $(E,\fm,\hh)$ with root system $R$, and introduce the notion of a ``cover" for a closed root subsystem $R'$ of $R$. A {\it cover} is a subspace $\hh' \subseteq \hh$ 
that preserves the non-degeneracy of the form and contains every coroot of $R'$ (Definition \ref{def:1}). In this case, the subalgebra $E_{R',\hh'}$ generated by $\hh'$ 
and the root spaces corresponding to non-zero roots of $R'$ forms an extended affine subalgebra 
called the \emph{Lie cover} of $(R',\hh')$ (Proposition \ref{pro:khan2}). 
We establish the existence of a canonical cover $\hh_R$ for the entire root system $R$, yielding a minimal-dimensional Lie cover $E_{R,\hh_R}$ 
that has the same core of $E$, and is tame if $E$ is tame (Corollary \ref{cor:prosiz1}). 
 Despite this, in general, the tameness does not inherits to a local Lie cover $E_{R',\hh'}$. To address this,
we introduce a modified construction by intersecting isotropic root spaces with the core, 
resulting in a triple $(\hat{E},\fm,\hh_R)$. 
Proposition \ref{pro:siz2} shows that this triple forms a tame EALA with root system $R$ 
if and only if the center of the core is contained in the Cartan subalgebra, and Corollary \ref{cor:6} establishes that tameness is guaranteed when, for example, these 
intersections are one-dimensional.

Section \ref{sec:affinization} develops the theory of affine localization. We begin by reviewing the affinization process from \cite{ABP02} (with further elaboration in \cite{AHY13}), which extends the theory behind the construction of twisted affine Lie algebras. Building upon these affinization techniques, we establish a crucial structural result concerning the intersection of isotropic root spaces with the core (Proposition \ref{pro:2}). This then leads to the central achievement of this paper: Given any affine root subsystem of $R$, there exists an associated local affine Lie subalgebra (Theorem \ref{thm:6} and Corollaries \ref{cor:10}-\ref{cor:12}). The section concludes with a complementary result concerning ``finite localization''.

Section \ref{examples} applies our systematic construction of affine Lie subalgebras to specific examples of extended affine Lie algebras, with particular attention to the detailed structural analysis of non-zero isotropic root spaces.

Section \ref{sec:filter} presents an application of the localization results, especially the existence of ascending filtrations for extended affine Lie algebras. Roughly speaking, starting from an EALA $E$ (subject to a mild restriction on types $A_1$ and $B_\ell$) whose root system $R$ has nullity $\nu$, we construct a filtration $R_0 \subseteq R_1 \subseteq \cdots \subseteq R_\nu = R$ of the root system $R$, where each $R_k$ has the same type and rank as $R$ but has nullity $k$. This filtration of the root system naturally gives rise to a corresponding filtration of Lie algebras $E_0 \subseteq E_1 \subseteq \cdots \subseteq E_\nu = E$ with each $E_k$ being an EALA whose root system is $R_k$ (Proposition \ref{pro:1}). A thorough investigation of these filtrations requires a separate study.

\markboth{S. Azam}{Chevalley Bases}

\section{\bf Preliminaries}\setcounter{equation}{0}\label{preliminaries}

 In this work, all vector spaces are considered over real or complex numbers. For a vector space $\v$ equipped with a symmetric bilinear form $\fm$, and a subset
$T$ of $\v$, we set
$$\begin{array}{l}
	T^\times:=\{\a\in \v\mid (\a,\a)\not=0\},\\
	T^0:=\v\setminus T^\times,\\
	 T^\perp:=\{\a\in\v\mid (\a,T)=\{0\}\}.
	 \end{array}$$
	  A subset $T\sub\v^\times$ is called {\it connected} if $T$ cannot be written in the form $T=T_1\cup T_2$ with
	 $(T_1,T_2)=\{0\}$ where $T_1\not=\emptyset$ and $T_2\not=\emptyset$.
	 
	 For $\a\in \v^\times$, we set 
$$\a^\vee:=\frac{2\a}{(\a,\a)}.$$
The dual space of a vector space $\u$ is denoted by $\u^\star$.

Given a subset $S$ of $\v$, we denote by $\la S\ra$ the additive subgroup generated by $S$.  

\subsection{Extended affine root systems}
We recall the definition and some basic properties of extended affine root systems from \cite{AABGP97}.
\begin{DEF}\label{def:azam1}
	An {\it extended affine root system} is a triple
	$(R,\fm,\v)$ where $\v$ is a finite dimensional real vector space, $\fm$ is a symmetric positive semi-definite bilinear form, and $R$ is a subset of $\v$ satisfying the following axioms:
	
	(R1) $0\in R$,
	
	(R2) $R=-R$,
	
	(R3) $R$ spans $\v$,
	
	(R4) $\a\in R^\times\Rightarrow 2\a\not\in R$,
	
	(R5) $R$ is discrete in $\v$,
	
	(R6) the {\it root string property}: for $\a\in R^\times$ and $\b\in R$, there exist non-negative integers $u,d$ such that
	$$(\b+\bbbz\a)\cap R=\{\b-d\a,\ldots,\b,\ldots,\b+u\a\}$$
	with $(\b,\a^\vee)=d-u$.
	
(R7) elements of $R^0$ are {\it non-isolated}, i.e., for $\sg\in R^0$ there exists $\a\in R^\times$ with $\a+\sg\in R$.


	(R8)  $R^\times$ is {\it connected}.
\end{DEF}

Let $(R,\fm,\v)$ be an extended affine root system. Define $\bar{\v} := \v / \v^0$, and let $\bar{\;} : \v \to \bar{\v}$ be the canonical projection. Then the image $\bar{R}$ of $R$ under this map is an irreducible finite root system in $\bar{\v}$, with respect to the bilinear form induced by $\fm$ on $\bar{\v}$. The \textit{type} and \textit{rank} of $R$ are defined to be the type and rank of $\bar{R}$, respectively. The \textit{nullity} of $R$ is defined to be the dimension of $\v^0$. {\it Throughout this work, we assume that $R$ is of reduced type}.

\begin{parag}\label{pt1}
	Given an EARS $(R,\fm, \v)$, it follows that $R$ contains a subset $\dot{R} \subseteq \v$ such that $R \subseteq \dot{R} + \v^0$, where $\dot{R}$ is a finite root system in $\dot{\v}:=\span_{\bbbr} \dot{R}$, and is isomorphic to $\bar{R}$ under the canonical projection. More specifically, we have:
	$$
	R = (S + S) \cup (\rds + S) \cup (\rdl + L), \quad \text{where } R^0 = S + S,
	$$
	for some subsets $S, L \subseteq \v^0$, called {\it semilattices}. Here $\rds$ and $\rdl$ are the set of short and long roots of $\rd$, respectively. If $\rd$ is simply laced, then every root is regarded as short (i.e., $\rd^\times = \rds$), and the term $\rdl + L$ is considered to be empty. If $\rdl\not=\emptyset$, the semilattices $S$ and $L$ interact as
	\begin{equation}\label{eq:inter1}
		\la L\ra+S= S\andd 2\la S\ra+L= L.
		\end{equation} 
		For more details about semilattices we refer the reader to \cite[Chapter II.\S 1]{AABGP97}
		and \cite{Az97}.

\begin{DEF}\label{def:n1}
	Let $(R,\fm, \v)$ be an extended affine root system.
	
	(i)
	We call a subset $R'$ of $R$ a {\it subsystem} of $R$ if the triple $(R',\fm',\v')$ is an extended affine root system, where
	$\v'=\span_{\bbbr}R'$ and $\fm'$ is the restriction to $\v'$ of $\fm$. In other words, an extended affine root system $(R',\fm',\v')$ is called a subsystem of $(R,\fm,\v)$ if $R'\sub R$ and $\fm'=\fm_{|_{\v'}}.$ 
	To indicate that $(R',\fm',\v')$ is a subsystem of $(R,\fm,\v)$, we simply write $(R',\fm',\v')\sub(R, \fm,\v)$.
	
	(ii) We call a subsystem $R'$ of $R$ {\it closed} if $\a,\b\in R'$, $\a+\b\in R$ implies $\a+\b\in R'$. Symbolically this is expressed as $(R'+R')\cap R\sub R'$. We further define $R'$ to be a {\it real-closed root subsystem}
	if for any $\a,\b\in\ R'$ with $\a+\b\in R$, the following implication holds:
	If at least one of 
	$\a$ or $\b$ is non-isotropic then $\a+\b\in R'$. In other words, if 
	$R'$ is closed under root addition whenever the sum is a root and at least one summand is non-isotropic. Symbolically, this can be stated as $(R'+{R'}^\times)\cap R\sub R'$.
\end{DEF}

We note that a subsystem of an irreducible root system, in general, is not closed even when it is finite and irreducible. A classical example is the $G_2$ root system. Its short roots form an irreducible $A_2$ subsystem. However, this subsystem is not closed because the sum of two short roots can be a long root of $G_2$. The study of ``closed" and ``real-closed" subsystems requires independent consideration; see for example \cites{Dyn52, BHV23, KV21}. 

\begin{parag}\label{parag:n1}
	There are some canonical ways of constructing a subsystem of an EARS $R$. Here, we describe one that is relevant to our purpose. Let $T$ be a non-empty connected subset of $R^\times$. Set
\begin{equation}\label{eq:n1}
	\begin{array}{l}
		\v_T:=\span_{\bbbr} T,\\
\tilde{R}_T^\times:=\la T\ra\cap R^\times,\\
\tilde{R}_T^0:=\v^0\cap (\tilde{R}_T^\times-\tilde{R}_T^\times),\\
\tilde{R}_T:=\tilde{R}_T^\times\cup \tilde{R}_T^0,\\
\fm_T:=\fm_{|_{\v_T}}.
\end{array}
		\end{equation}
	It is straightforward to check that $(\v_T, \fm_T, {\tilde R}_T)$ is a real-closed subsystem of $R$. Moreover, $({\tilde R}_T,\fm_T,\v_T)$ is the maximal subsystem of $R$ containing $T$.
	
	We observe that if $\tilde{R}_T^0$ forms a lattice, then Definition \ref{def:n1} ensures that the properties of being closed and real-closed are equivalent. In particular, since for any extended affine root system of nullity $\leq 2$, the isotropic roots form a lattice, $\tilde{R}_T$ is closed if it has nullity $\leq 2$.
\end{parag}
%
%
%
%
\end{parag}

\begin{exa}\label{exa:-1}
	Let $(R,\fm,\v)$ be an extended affine root system and $\dot R$ and $\dot\v$ be as in \pref{pt1}. 
	
	(i) Let $0\not=\d\in R^0$ and set $T=({\dot R}^\times+\bbbz\d)$. It follows from \pref{pt1} that
	$T\cap R^\times\not=\emptyset$. As indicated in \pref{parag:n1}, 
	$(\tilde{R}_T,\fm_T,\v_T)$ is an extended affine root system. 
	We have
	$$\tilde{R}_T^\times =\la T\ra\cap R^\times= (\la\rd\ra\oplus\bbbz\d)\cap R^\times.$$
	Since $\tilde{R}_T^0\sub\bbbz\d$,  $\tilde R$ has nullity $1$ and so is a closed affine root subsystem of $R$, with the same type as $R$. We denote this affine root subsystem by $R_{\dot R,\d}$. 
	
	(ii) Let $\a\in R^\times$. Then $\a=\dot\a+\d$ for some $\dot\a\in\rd^\times$ and some $\d\in R^0$. Set
	$T=\{\dot\a,\a\}$. Then by \pref{parag:n1}, $\tilde{R}_T$ is a closed affine subsystem of $R$ containing $\a$.
	We denote this subsystem by $R_{\a,\d}$.
	\end{exa}

\subsection{\bf Extended affine Lie algebras}\setcounter{equation}{0}\label{subsec:1}
We recall the definition of an extended affine Lie algebra and  highlight some of its fundamental properties. Unless stated otherwise, the results discussed here are drawn from \cite[Chapter I]{AABGP97}. All algebras in this work are considered over filed $\bbbk$ of complex numbers.

\begin{DEF}\label{def:ea1}
	An {\it extended affine Lie algebra} is  a triple $(E,\fm,\hh)$ where $E$ is a Lie algebra, $\hh$ is a subalgebra of $E$ and $\fm$ is a bilinear form on $E$  satisfying the following 5 axioms:
	
	(EA1) The form $\fm$ on $E$ is symmetric, non-degenerate and invariant.
	
	(EA2)  $\hh$ is a finite dimensional splitting Cartan subalgebra of $E$. This means that $E=\sum_{\a\in\hh^\star}E_\a$ where $E_\a=\{x\in E\mid [h,x]=\a(h)x\hbox{ for all }h\in\hh\}$ and $E_0=\hh$. 
	
	
	Let $R$ be the set of roots of $E$, namely $R=\{\a\in\hh^\star\mid E_\a\neq\{0\}\}$.
	It follows from (EA1)-(EA2) that the form $\fm$ restricted to $\hh$ is non-degenerate and so it can be transferred
	to $\hh^\star$ by $(\a,\b):=(t_\a,t_\b)$ where $t_\a\in\hh$
	is the unique element satisfying $\a(h)=(h,t_\a)$, $h\in\hh$.
	Then $R=R^0\uplus R^\times$ is regarded as the decomposition of roots into
	{\it isotropic} and {\it non-isotropic} roots, respectively. Let $\v:=\hbox{Span}_\bbbr R$
	and $\v^0:=\hbox{Span}_{\bbbr}R^0$.
	
	(EA3) For $\a\in R^\times$, $\ad x$ is locally nilpotent for $x\in E_\a$.
	
	(EA4) $R$ is discrete in $\hh^\star$.
		
	(EA5) $R$ is irreducible. That is, $R^\times$ is connected and elements of $R^0$ are non-isolated (see axioms (R7) and (R8) of Definition \ref{def:azam1}).
	
%
%
	\end{DEF}
	
Let $(E,\fm,\hh)$ be an extended affine Lie algebra with root system $R$. It follows that $(R,\fm,\v)$ is an extended affine root system in the sense of Definition \ref{def:azam1}. 
We have 
\begin{equation}\label{eq:khat0}
[E_\a,E_{-\a}]=\bbbk t_\a,\quad(\a\in R),
\end{equation}
\begin{equation}\label{eq:khat1}
	(E_\a,E_\b)=\{0\}\hbox{ unless }\a+\b=0,\quad(\a,\b\in R).
	\end{equation} 
	Also from \cite[Remark 1.5]{Az06}, we have
	\begin{equation}\label{eq:khat4}
		[E_\a, E_\b]\not=\{0\},\qquad(\a\in R^\times,\b,\a+\b\in R).
		\end{equation}

The {\it core} $E_c$ of $E$ is by definition the subalgebra of $E$ generated by the non-isotropic root spaces. It follows that the centralizer $C_E(E_c)$ of  $E_c$ in $E$ coincides with the orthogonal complement $E_c^\perp$ of $E_c$ in $E$. Therefore 
\begin{equation}\label{eq:9}
	Z(E_c)\sub C_E(E_c)=E_c^\perp,
	\end{equation}
	and
	\begin{equation}\label{eq:9-1}
		x\in E_c^\perp\Longleftrightarrow (x,E_\a)=\{0\}\hbox{ for all }\a\in R^\times.
		\end{equation}
	
	Note that if $x\in E_c^\perp$, then  by (\ref{eq:khat1}) each homogeneous component of $x$ is contained in $E_c^\perp$. It then follows that
	\begin{equation}\label{eq:omid1}
		E_c^\perp\sub\sum_{\sg\in R^0}{E_\sg}.
		\end{equation}
	 The extended affine Lie algebra $E$ is called {\it tame} if $E_c^\perp=Z(E_c)$ (or equivalently $E_c^\perp\sub Z(E_c)$).

The following lemma gives a helpful characterization of tameness.

\begin{lem}\label{lem:tame3}
	The extended affine Lie algebra $E$ is tame if and only if for each $\sg\in R^0$,
	$$E_\sg\cap E_c^\perp=Z(E_c)\cap \sum_{\a\in R^\times}[E_{\a+\sg},E_{-\a}].$$
\end{lem}

\proof We have
\begin{eqnarray*}
	E\hbox{ is tame}&\Longleftrightarrow& E_c^\perp=Z(E_c)\\
	&\Longleftrightarrow& E_c^\perp\cap E_\a=Z(E_c)\cap E_\a,\;\forall\a\in R\\
	&\Longleftrightarrow& E_c^\perp\cap E_\sg=Z(E_c)\cap E_\sg,\;\forall\sg\in R^0\\
	&\Longleftrightarrow&E_c^\perp\cap E_\sg=Z(E_c)\cap \sum_{\a\in R^\times}[E_{\a+\sg}, E_{-\a}],\;\forall\sg\in R^0.\\
\end{eqnarray*}
\qed

\begin{lem}\label{eq:naz4}
For $\a\in R^\times$ and $\sg\in R^0$ with $\a+\sg\in R$, we have
$$([E_{\a+\sg},E_{-\a}],[E_{-\a-\sg},E_{\a}])\not=\{0\}.$$
 \end{lem}
 
 \proof Take $x_{\pm\a}\in E_{\pm\a}$ and 
 $x_{\pm(\a+\sg)}\in E_{\pm(\a+\sg)}$ with $[x_\a,x_{-\a}]=t_\a$ and $[x_{\a+\sg},x_{-\a-\sg}]=t_{\a+\sg}.$ From the Jacobi identity, the fact that $-2\a-\sg$ is not a root and invariance of the form, we have 
\begin{eqnarray*}
 ([x_{\a+\sg},x_{-\a}],[x_{-\a-\sg},x_{\a}])
 &=& (x_{\a+\sg},[x_{-\a},[x_{-\a-\sg},x_{\a}]])\\
 &=&
 -(x_{\a+\sg}, [x_{-\a-\sg},[x_\a,x_{-\a}]])\\
 &=&
 -(x_{\a+\sg},[x_{-\a-\sg},t_\a])\\
 &=&
 -([x_{\a+\sg},x_{-\a-\sg}],t_\a)\\
 &=&-(t_{\a+\sg},t_\a)=-(\a,\a)\not=0.
 \end{eqnarray*}
 \qed

\section{A Construction of Extended Affine Lie Algebras}
\setcounter{equation}{0}\label{sec:tame6}

In this section, we provide a brief overview of a general approach to constructing extended affine Lie algebras, following the work of Erhard Neher \cite{Neh04}. We assume the reader is already familiar with the concept of Lie tori (see \cite{Yos06}).  

Let $\Lambda$ be a lattice of rank $\nu$, and let $\mathfrak{g}$ be a centerless $\Lambda$-Lie torus, meaning a Lie torus with a trivial center. It is known that $\fg$ admits a symmetric invariant non-degenerate bilinear from $\fm_\fg$. Consider the \textit{centroid} $\cent_\mathbb{K}(\mathfrak{g})$ of $\mathfrak{g}$, which is $\Lam$-graded with
$\dim_\mathbb{K} \cent_\mathbb{K}(\mathfrak{g})^\lambda \leq 1.$
Defining  
$\Gamma := \{ \lambda \in \Lambda \mid \cent_\mathbb{K}(\mathfrak{g})^\lambda \neq 0 \},
$
we can express the centroid as  
$
\cent_\mathbb{K}(\mathfrak{g}) = \bigoplus_{\mu \in \Gamma} \mathbb{K} \chi^\mu,
$
where each $\chi^\mu$ acts as an endomorphism of degree $\mu$, satisfying $\chi^\mu \chi^\nu = \chi^{\mu+\nu}$.  

For any $\theta \in \text{Hom}_\mathbb{Z}(\Lambda, \mathbb{K})$, the \textit{degree derivation} $\partial_\theta$ on $\mathfrak{g}$ is defined by  
$
\partial_\theta(x^\lambda) = \theta(\lambda) x^\lambda,$ for  $x^\lambda \in \mathfrak{g}^\lambda.
$
Denoting by $\mathcal{D}$ the set of all degree derivations, the set of \textit{centroidal derivations} of $\mathfrak{g}$ is defined by  
$
\mathrm{CDer}_\mathbb{K}(\mathfrak{g}) := \cent_\mathbb{K}(\mathfrak{g}) \mathcal{D} = \bigoplus_{\mu \in \Gamma} \chi^\mu \mathcal{D}.
$
This forms a $\Gamma$-graded subalgebra of the derivation algebra $\mathrm{Der}_\mathbb{K}(\mathfrak{g})$, with the bracket  
\begin{equation}\label{eq7}
	[\chi^\mu \partial_\theta, \chi^\nu \partial_\psi] = \chi^{\mu+\nu} (\theta(\nu) \partial_\psi - \psi(\mu) \partial_\theta).
\end{equation}

Next, the algebra of \textit{skew centroidal derivations} of $\mathfrak{g}$ is given by
$$
\mathrm{SCDer}_\mathbb{K}(\mathfrak{g}) = \bigoplus_{\mu \in \Gamma} \mathrm{SCDer}_\mathbb{K}(\mathfrak{g})^\mu = \bigoplus_{\mu \in \Gamma} \chi^\mu \{ \partial_\theta \in \mathcal{D} \mid \theta(\mu) = 0 \}.
$$
We note that $\mathrm{SCDer}_\mathbb{K}(\mathfrak{g})^0 = \mathcal{D}$.  

For a graded subalgebra $D = \sum_{\mu \in \Gamma} D^\mu$ of $\mathrm{SCDer}_\mathbb{K}(\mathfrak{g})$, its \textit{graded dual} is defined as  
\[
D^{gr^\star} = \sum_{\mu \in \Gamma} (D^\mu)^\star, \quad \text{where } (D^{gr^\star})^\mu = (D^{-\mu})^\star.
\]
This dual space is treated as a $D$-module via  
$
(d.\varphi)(d') = \varphi([d', d]),$ for  $d', d \in D$, $\varphi \in D^{gr^\star}.
$
Here, $\varphi \in (D^\mu)^\star$ is viewed as a linear map on $D$ with $\varphi|_{D^\nu} = 0$ for $\nu \neq \mu$.  

A $\Gamma$-graded subalgebra $D = \bigoplus_{\mu \in \Gamma} D^\mu$ of $\mathrm{SCDer}_\mathbb{K}(\mathfrak{g})$ is called \textit{permissible} if the canonical evaluation map  
$\text{ev}: \Lambda \to (D^0)^\star$, defined by $\text{ev}(\lambda)(\partial_\theta) = \theta(\lambda),
$
is injective and has a discrete image.  

Finally, let $\kappa$ be a bilinear map satisfying  
$$
\begin{array}{c}
	\kappa(d, d) = 0, \quad \sum_{(i,j,k) \circlearrowleft} \kappa([d_i, d_j], d_k) = \sum_{(i,j,k) \circlearrowleft} d_i \cdot \kappa(d_j, d_k), \vspace{2mm} \\
	\kappa(D^{\mu_1}, D^{\mu_2}) \subseteq (D^{-\mu_1-\mu_2})^\star, \quad \kappa(d_1, d_2)(d_3) = \kappa(d_2, d_3)(d_1), \vspace{2mm} \\
	\kappa(D^0, D) = 0,
\end{array}
$$
for $ d, d_1, d_2, d_3 \in D $. The map $\kappa$ is referred to as an \textit{affine cocycle} on $D$.  

Now, assuming $\mathfrak{g}$ is a centerless Lie torus, $D$ is a permissible subalgebra of $\mathrm{SCDer}(\mathfrak{g})$, and $\kappa$ is an affine cocycle on $D$, then the vector space  
\begin{equation}\label{rok1}
	E = E(\mathfrak{g}, D, \kappa) := \mathfrak{g} \oplus D^{gr^\star} \oplus D,
\end{equation}
forms a Lie algebra with the bracket
\begin{eqnarray*}
	[x_1+c_1+d_1,x_2+c_2+d_2]&=&([x_1,x_2]_\fg+d_1(x_2)-d_2(x_1))\\
	&+&(c_D(x_1,x_2)+d_1.c_2-d_2.c_1+\kappa(d_1,d_2))\\
	&+&[d_1,d_2]
\end{eqnarray*}
for $x_1,x_2\in \fg,c_1,c_2\in D^{gr*},d_1,d_2\in D$. Here $[\:,\:]_\fg$ denotes the Lie bracket on $\fg$, $[d_1,d_2]=d_1d_2-d_2d_1$, and $c_D:\fg\times\fg\rightarrow D^{gr*}$ is defined by
$$c_D(x,y)(d)=(d(x)|y)\;\text{for all}\;x,y\in \fg,d\in D.$$
Finally, the bilinear form $\fm$ on $E$ defined by
\begin{equation}\label{ira1}
	(x_1+c_1+d_1,x_2+c_2+d_2)=(x_1,x_2)_\fg+c_1(d_2)+c_2(d_1)
\end{equation}
is symmetric, invariant and non-degenerate.

This construction yields a \textit{tame extended affine Lie algebra} $(E,\fm, \hh)$, where  
$
\hh= \mathfrak{g}^0_0 \oplus (D^0)^\star \oplus D^0.
$
To indicate the dependence of $E$ on $\fg$, $D$ and $\kappa$, we write
$E=E(\fg,D,\kappa)$. Note that
$E_c=\fg\oplus {D^{gr}}^\star.
$
From the way the bracket is defined on $E$, we see that
\begin{equation}\label{eq:tame4}
	Z(E_c)={D^{gr}}^\star,
\end{equation}
and
	\begin{equation}\label{eq:tame8}
	Z(E_c)\sub\hh\Longleftrightarrow D=D^0\Longleftrightarrow Z(E_c)={D^0}^\star.
\end{equation}	
\begin{rem}\label{eq:tame7}
Most interesting examples of extended affine Lie algebras in the literature have $D=D^0$ and $\kappa=0$, as one can see, for example, in \cite[Chapter III]{AABGP97}, \cite{BGKN95}, \cite{BGK96} and \cite{H-KT90}.
\end{rem}

\section{\bf Subalgebras and Localization}\setcounter{equation}{0}\label{sec:local}
Let $(E,\fm,\hh)$ be an extended affine Lie algebra with root system $R$. We recall from Section \ref{subsec:1} that the form on $\hh$ induces  a form on
$\hh^\star$ by $(\a,\b)=(t_\a,t_\b)$, $\a,\b\in R$.	It then follows that $(R,\fm,\v)$ is an extended affine root system where
$\v=\span_\bbbr R$ and $\fm$ is the form on $\hh^\star$ restricted to $\v$.

\subsection{Local Lie covers}
Assume that $(E,\fm,\hh)$ is an EALA with root system $R$.

\begin{DEF}\label{def:1} Let $R'$ be a closed subsystem of $R$.
	
(i) We call a subspace $\hh'$ of $\hh$ a {\it (local) cover} for $R'$ if $t_\a\in \hh'$ for all $\a\in {R'}^\times$, and the form on $\hh$ restricted to $\hh'$ is non-degenerate.

(ii) For a cover $\hh'$ of $R'$, we set
$$E_{R',\hh'}:=\hh'\oplus\sum_{\a\in R'\setminus\{0\}}E_\a,$$ and consider it as a vector space equipped with the form inherited from $\fm$. We call $E_{R',\hh'}$ the {\it (local) Lie cover} associated with  $(R',\hh')$ in $E$. The term Lie cover will be justified in Proposition \ref{pro:khan2}.
\end{DEF}


\begin{parag}\label{const1} Let $R'$ be a subsystem of $R$, and $\hh'$ be a cover for $R'$. 
Since each isotropic root in $R'$ is attached to some non-isotroic root in $R'$ (see (R8)), we have $t_\a\in\hh'$ for all $\a\in R'$.
Also since the forms on $\hh'$ and $\hh$ are non-degenerate, we have $\hh=\hh'\oplus{\hh'}^\perp.$
				Now since $(\hh',\hh'^\perp)=\{0\}$, and for $\a\in R'$, $t_\a\in\hh'$, we may identify $\a$ with its restriction to $\hh'$ and consider $\a$ as an element of ${\hh'}^\star$.  
				\begin{equation}\label{eq:identify}
				\a\in R'\Longrightarrow	t_\a\in\hh'\andd\a\equiv\a_{|_{\hh'}}.
				\end{equation}
				This gives
			\begin{equation}\label{eq:identify2}
			{R'}^\times=R'\cap R^\times\andd {R'}^0=R'\cap R^0.
			\end{equation}
			The facts (\ref{eq:identify})-(\ref{eq:identify2}) will be used frequently in the sequel without further reference.
\end{parag}

\begin{pro}\label{pro:khan2}
	Let $R'$ be a closed subsystem of $R$, and $\hh'$ a cover of $R'$. Then the Lie cover $E_{R',\hh'}$ is an extended affine Lie algebra with Cartan subalgebra $\hh'$ and root system $R'$. In particular, $E_{R',\hh}$ is an extended affine Lie algebra with root system $R'$.
\end{pro}

\proof Set $E'=E_{R',\hh'}$. Since $\hh'$ is abelian and $R'$ is a closed subsystem of $R$, $E'$ is a Lie subalgebra of $E$.
Since the form on $E$ is non-degenerate, it follows from (\ref{eq:khat1}) that the form on $E'$ is also non-degenerate, thus (EA1) holds for $E'$. 

Next, since $E'$ is an $\hh$-submodule of $E$, we have $E'=\sum_{\a\in R}^\oplus E'_\a=\sum_{\a\in R'}E'_\a$ where
$E'_\a=E'\cap E_\a$. Now, we show that for $\a\in R'$, 
$$E_\a\cap E'=\{x\in E'\mid[h,x]=\a(h)x\hbox{ for all }h\in\hh'\},$$
where
the inclusion $``\sub"$ is clear. To see the reverse inclusion, let $x\in E'$ and $[h',x]=\a(h')x$ for all $h'\in\hh'$.
We have $x=x_{\a_1}+\cdots+ x_{\a_k}$, for some distinct $\a_i\in R'$, $x_{\a_i}\in E'_{\a_i}$. Now let $h\in\hh$. Then, $h=h'+h''$ where $h'\in\hh'$ and $h''\in{\hh'}^\perp$. Since $\a({\hh'}^\perp)=\a_i({\hh'}^\perp)=\{0\}$ for all $i$, we get 
\begin{eqnarray*}\a(h)x-[h,x]&=&
\a(h')x-\sum_{i=1}^k\a_i(h'+h'')x_{\a_i}\\
	&=&[h',x]-\sum_{i=1}^k\a_i(h')x_{\a_i}\\
	&=&[h',x]-\sum_{i=1}^k[h',x_{\a_i}]=0.
	\end{eqnarray*}
This gives $x\in E_\a\cap E'$. It is easy to see that $\hh'$ is self-centralizing in $E'$. Thus (EA2) holds for 
$(E',\fm,\hh')$ and $R'$ is the root system of $E'$ with respect to $\hh'$.
 Axioms (EA3)-(EA5) trivially hold since
$R'$ is a subsystem of $R$.
\qed

\subsection{A canonical Lie cover}
Let $(E,\fm,\hh)$ be an EALA with root system $R$. The following notation will be used frequently in the sequel:

\noindent{\bf Notation.}\label{not:1}
 For a subsystem $R'$ of $R$, we set 
$$\hh_c(R'):=\sum_{\a\in {R'}^\times}\bbbk t_\a\andd \hh^0(R'):=\sum_{\sg\in {R'}^0}\bbbk t_\a.$$	
For notional convenience, when $R'=R$ we abbreviate  $\hh_c(R)$ and $\hh^0(R)$ by $\hh_c$ and $\hh^0$, respectively.

\begin{const}\label{const:1} We construct a canonical Lie cover for the extended affine Lie algebra $(E,\fm,\hh)$ with a minimal property as will be illustrated in Corollary \ref{cor:prosiz1} below. 

	We have
\begin{equation}\label{eq:khat2}
	\hh\cap E_c=E_0\cap E_c=\sum_{\a\in R^\times}[E_\a,E_{-\a}]=\sum_{\a\in R^\times}\bbbk t_\a=\hh_c.
	\end{equation}

Let $\rd$ be as in \pref{pt1}. We set 
$$\dot\hh:=\sum_{\dot\a\in\rd}\bbbk t_{\dot\a}.$$
Since each $\a\in R$ is of the form $\dot\a+\sg$ for some
$\dot\a\in\rd$ and some $\sg\in R^0$, we have
$$\hh_c=\dot\hh\oplus\hh^0.$$
Since the form on both $\dot\hh$ and $\hh$ is non-degenerate, we have $\hh=\dot\hh\oplus{\dot\hh}^\perp,$ where the form on ${\dot\hh}^\perp$ being also non-degenerate, and $\hh^0\sub\dot\hh^\perp.$ By Witt's extension theorem \cite[Theorem 11.11]{Rom05}, $\dot\hh^\perp$ contains a subspace 
$\hat\hh^0$ such that 
\begin{equation}\label{eq:khat3}
	\begin{array}{c}
		\dim\hat\hh^0=\dim\hh^0,\\
		\hbox{the form restricted to $\hh^0\oplus\hat\hh^0$ is non-degenerate},\\
		(\hat\hh^0,\hat\hh^0)=\{0\}.
	\end{array}
\end{equation}

We set 
$$\hh_R:=\hh_c\oplus\hat\hh^0=\dot\hh\oplus\hh^0\oplus{\hat\hh}^0.$$
Clearly the form on $\hh_R$ is non-degenerate. 
Therefore,  
\begin{equation}\label{eq:tame2}
	\hh=\hh_R\oplus\hh_R^\perp.
\end{equation}
To simplify the notation, we denote by $E_R$ the Lie cover $E_{R,\hh_R}$ associated to the pair $(R,\hh_R)$.
\end{const}

\begin{cor}\label{cor:prosiz1}
The Lie cover $E_R$ is an extended affine Lie algebra with Cartan subalgebra $\hh_R$, root system $R$, and core $(E_R)_c=E_c$. It has minimal dimension among all Lie covers containing the core. Moreover, if $E$ is tame then so is $E_R$.
\end{cor}

\proof
 Since $\hh_R$ is a cover for $R$, the first part of statement is an immediate consequence of Definition \ref{def:1} and Proposition \ref{pro:khan2}. Now the Cartan subalgebra $\hh'$ of any Lie cover containing the core contains $\hh_c$ and so the non-degeneracy of the form on $\hh'$ forces $\dim \hh'\geq\dim \hh_R$.
 
 If $E$ is tame and $x\in (E_R)_c^\perp$, then $x\in E_R$ and
 $(x,E_c)=\{0\}$, so $x\in E_c^\perp=Z(E_c)=Z((E_R)_c).$\qed

\begin{rem}\label{rem:m1} 
  As previously noted, all extended affine Lie algebras and root systems in this work are assumed to be reduced. Nevertheless, Proposition \ref{pro:khan2} and Corollary \ref{cor:prosiz1} remain valid even for non-reduced extended affine Lie algebras, since their proofs do not rely on the reduceness condition.
 	\end{rem}

\subsection{Core-based localization} Let $(E,\fm,\hh)$ be an EALA with root system $R$. 
 Despite Corollary \ref{cor:prosiz1}, in general the tameness does not inherits to a local Lie cover $R_{R',\hh'}$. To address this, we consider another localization here.

We set
\begin{equation}\label{eq:rok2}
\hat{E}=\hh_R+E_c=(\hh_c\oplus\hat\hh^0)+ E_c=\hat\hh^0\oplus E_c.
\end{equation}
 
\begin{lem}\label{lem:1}
	If $E$ is tame then $\hh_R^\perp=\{0\}$ and $\hat E=\hh+E_c$.
	\end{lem}
	
	\proof If $E$ is tame then $E_c^\perp\sub E_c$ and so
$\hh_R^\perp\sub\hh\cap E_c^\perp\sub\hh\cap E_c= \hh_c.$ Thus $\hh_R^\perp=\{0\}$.
	Hence $\hh=\hh_R$ and $\hat{E}=\hh+E_c$.\qed

For further study of $\hat{E}$, we set 
	$$\hat{E}^{{iso}}:=(\sum_{0\not=\sg\in R^0}E_\sg)\cap E_c=\sum_{0\not=\sg\in R^0}\sum_{\a\in R^\times}[E_{\a+\sg},E_{-\a}].$$
	Then
	$$\hat{E}=\hh_R\oplus\hat{E}^{{iso}}\oplus\sum_{\a\in R^\times}E_\a\andd \hat{E}_c=\hh_c\oplus\hat{E}^{{iso}}\oplus\sum_{\a\in R^\times}E_\a.$$

\begin{lem}\label{lem:5} 
The form restricted to $\hat{E}$ is non-degenerate if and only if the form on $\hat{E}^{{iso}}$ is non-degenerate if and only if $Z(E_c)\sub\hh$. 
\end{lem}

\proof
The first if and only if is clear from (EA1), (\ref{eq:khat1}) and the fact that $\fm$ is non-degenerate on
$\hh_R$ and $\sum_{\a\in R^\times} E_\a$.

Next, the form on ${\hat E}^{iso}$ is non-degenerate if and only if the radical 
$K_\sg$ of the form restricted to $(E_\sg\cap E_c)\oplus (E_{-\sg}\cap E_c)$ is trivial for each 
$\sg\in R^0\setminus\{0\}$. Now, as for a fixed 
$\sg\in R^0\setminus\{0\}$, we have $(K_\sg,E_c)=\{0\}$, and $Z(E_c)=E_c\cap E_c^\perp=(\hh_c\oplus{\hat E}^{iso})\cap E_c^\perp$, thus
the form on ${\hat E}^{iso}$ is degenerate if and only if 
$$\begin{array}{c}
K_\sg\not=\{0\}\hbox{ for some }\sg\in R^0\setminus\{0\}\\
\Longleftrightarrow\\
\{0\}\not=K_\sg\sub E_c\cap E_c^\perp\cap {\hat E}^{iso}\hbox{ for some }\sg\in R^0\setminus\{0\}\\
\Longleftrightarrow\\
Z(E_c)\cap {\hat E}^{iso}\not=\{0\}.
\end{array}
$$
Thus the form on $\hat E$ is non-degenerate if and only if $Z(E_c)\cap {\hat E}^{iso}=\{0\}$ if and only if
$Z(E_c)\sub \hh$.

\begin{pro}\label{pro:siz2}
	Let $(E,\fm,\hh)$ be an extended affine Lie algebra with root system $R$. Then
	$(\hat{E},\fm, \hh_R)$ is a tame extended affine Lie algebra (with root system $R$) if and only if $Z(E_c)\sub\hh$.
\end{pro}

\proof If $Z(E_c)$ is not contained in $\hh$, then $\fm$ is degenerate by Lemma \ref{lem:5}, and so $(\hat E,\fm,\hh_R)$ is not an extended affine Lie algebra.

Conversely, assume that $Z(E_c)\sub\hh$. By Lemma \ref{lem:5}, the form on  $\hat E$ is non-degenerate ensuring that (EA1) holds for $\hat E$.
 
 Next, since isotropic roots are non-isolated, for each 
$\sg\in R^0$ there exists 
$\a\in R^\times$ with $\a+\sg\in R^\times$. Then by (\ref{eq:khat4}),
\begin{equation}\label{eq:siz1}
	\{0\}\not=[E_{\a+\sg},E_{-\a}]\sub E_\sg\cap E_c,\qquad(\sg\in R^0).
\end{equation}

By Corollary \ref{cor:prosiz1}, $E_{R}$ is an extended affine Lie algebra with root system $R$. Since $\hat{E}$ is an $\hh_{R}$-submodule of $E_{R}$, it inherits a weight space decomposition
$\hat{E}=\sum_{\a\in R}\hat{E}_\a$ with 
$$\hat{E}_\a=(E_{R})_\a\cap\hat{E}=\{x\in \hat{E}\mid [h,x]=\a(h)x\hbox{ for all }h\in\hh_{R}\}.$$
In particular $\hat E_0=(E_R)_0\cap \hat E=\hh_R\cap\hat E=\hh_R$.
By (\ref{eq:siz1}), the set of weights of $\hat{E}$ is also $R$. Thus (EA2) also holds. Axioms (EA3)-(EA5) trivially hold. Thus $(\hat{E},\fm,\hh_R)$ is an extended affine Lie algebra. 

Finally, we show that $\hat E$ is tame. Let $x\in \hat{E}_c^\perp.$ Then
$x=\sum_{\sg\in R^0}x_\sg$, where $x_\sg\in \hat{E}_c^\perp\cap \hat{E}_\sg$, $\sg\in R^0$. Since the form on $\hat{E}^{{iso}}$ is non-degenerate,
we conclude that $x=x_0$. Since $x_0\in \hh_R$, we have $x_0=\dot x+x^0+\hat{x}^0$, where $\dot x\in\dot\hh$, $x^0\in\hh^0$ and $\hat{x}^0\in\hat{\hh}^0$. Now since $(x_0,\dot\hh\oplus\hh^0)=\{0\}$, we conclude that $x=x_0=x^0\in\hh^0\sub\hat{E}_c$.
\qed

\begin{cor}\label{cor:6}
	Let $(E,\fm,\hh)$ be an extended affine Lie algebra with root system $R$. If for every nonzero $\sigma\in R^0\setminus\{0\}$ the spaces $E_{\pm\sigma}\cap E_c$ are $1$-dimensional, then:
	\begin{enumerate}[(i)]
		\item The form $\fm$ restricted to $(E_{\sigma}\cap E_c)\oplus (E_{-\sigma}\cap E_c)$ is non-degenerate, for each $\sg\in R^0\setminus\{0\}$.
		\item The triple $(\hat{E},\fm,\hh_R)$ is a tame extended affine Lie algebra.
	\end{enumerate}
\end{cor}

\proof
(i)	Let $\sigma\in R^0\setminus\{0\}$. Since isotropic roots are non-isolated, there exists $\alpha\in R^\times$ such that $\alpha+\sigma\in R^\times$. When the spaces $E_{\pm\sigma}\cap E_c$ are one-dimensional, equation (\ref{eq:khat4}) gives:
	\begin{eqnarray*}
		E_{\sigma}\cap E_c = [E_{\alpha+\sigma},E_{-\alpha}]\andd
		E_{-\sigma}\cap E_c = [E_{-\alpha-\sigma},E_{\alpha}].
	\end{eqnarray*}
	Lemma \ref{eq:naz4} then guarantees that the form $\fm$ is non-degenerate on 
	\[ \hat{E}_\sigma\oplus\hat{E}_{-\sigma} = [E_{\alpha+\sigma},E_{-\alpha}]\oplus [E_{-\alpha-\sigma},E_{\alpha}]. \]
	
(ii)	If the spaces $E_{\pm\sigma}\cap E_c$ are $1$-dimensional for all $\sigma\in R^0\setminus\{0\}$, then
by(i), the form $\fm$ is non-degenerate on $\hat{E}$. Lemma \ref{lem:5} consequently gives $Z(\hat{E}_c)\subseteq \hh_R$, and Proposition \ref{pro:siz2} establishes that $(\hat{E},\fm,\hh_R)$ is tame.
	\qed

\begin{rem}\label{rem:7} 
In the literature, nearly all instances of EALAs exhibit the property $Z(E_c)\sub \hh$; for example those discussed in \cite{H-KT90}, \cite{Pol94} and \cite[Chapter III]{AABGP97}.
	\end{rem}

\section{Affine Localization}\setcounter{equation}{0}\label{sec:affinization}
While our primary focus in this section is on affine localization, we establish a general framework applicable to extended affine Lie algebras arising from the so-called affinization process, which extends the classical construction of twisted affine Lie algebras. Using this framework, we analyze the interaction of isotropic root spaces with the core (Proposition \ref{pro:2}). The main result of this section (theorem \ref{thm:6}) establishes the existence of local affine Lie subalgebras within an extended affine Lie algebra.
For a detailed study of the affinization process, we refer the reader to \cite{ABP02} and \cite{AHY13}.

Let $(\fg,\fm,\fh)$ be an extended affine Lie algebra with root system $R$. 
Let $\sg\in\Aut(\fg)$ be an automorphism of $\fg$. We denote by $\fg^\sg$ and $\fh^\sg$ the fixed points of $\fg$ and $\fh$ under $\sg$, respectively. Let $m$ be a positive integer and that the $\sg$ satisfies the following conditions:
\begin{equation}\label{eq:axiom}
	\begin{array}{l}
\sg^m=\id,\\
 \sg(\fh)=\fh,\\
 (\sg(x),\sg(y))=(x,y)$ for all $x,y\in\fg,\\
 C_{\fg^\sg}(\fh^\sg)=\fh^\sg.
\end{array}
\end{equation}
(Here $C_{\fg^\sg}(\fh^\sg)$ denotes the centralizer of $\fh^\sg$ in $\fg^\sg$.) 

Consider the group epimorphism $\bbbz\rightarrow\bbbz_m$ with $i\mapsto\bar i$. Let $\omega$ be a primitive $m$-th root of unity.
Then
$$
\fg=\sum_{\bar{i}\in\bbbz_m}\fg^{\bar i}\quad\hbox{ where }\quad\fg^{\bar i}:=\{x\in\fg\mid \sg(x)=\omega^ix\}.$$
Since $\sg$ preserves $\fh$, we have a similar decomposition for $\fh$.
Let we denote the transpose of $\sg$ on $\fh^\star$ by $\sg$ again. In other words,
$\sg(\a)(h)=\a(\sg^{-1}(h))$, $h\in\fh$, $\a\in\fh^\star$. To simplify the notation, and since distinction is clear from the context, we use the same symbol  $\pi_{\bar j}$ for all three projections below:
$$\fg\rightarrow\fg^{\bar j},\quad \fh\rightarrow\fh^{\bar j},\quad \fh^\star\rightarrow{\fh^\star}^{\bar j}.$$
		We also note that $\sg$ stabilizes the core $\fg_c$ of $\fg$ and so it stabilizes $\fh_c=\fh\cap\fg_c$. Therefore the projection $\pi_{\bar j} :\fh_c\rightarrow\fh_c^{\bar j}$ also makes sense.
		 For the sake of notation we denote $\pi_{\bar 0}$ by $\pi$, and $\pi_{\bar j}$ by $\pi_j$, $\bar j\not=\bar 0$. Then
		\begin{equation}\label{eq:pi1}
			\pi_j=\frac{1}{m}\sum_{i=0}^{m-1}\omega^{-ij}\sg^i.
		\end{equation} 
		 Here are some useful formulas concerning the projections $\pi_j$ which will be of use in the sequel, (see \cite[\S 3]{AHY13}). For $0\leq j,n\leq m-1$, and $x,y\in\fg$,
	 	\begin{equation}\label{eq:4}
	 	\begin{array}{c}
	 			\sg\circ\pi_j=\pi_j\circ\sg=\omega^j\pi_j,\vspace{2mm}\\  
	 	\;	[\pi_n(x),\pi_j(y)]=\pi_{n+j}[x,\pi_j(y)],\vspace{2mm}\\
	 	\;(\pi(x),\pi(y))=(x,\pi(y)).
	 		\end{array}
	 		\end{equation}

Next let $\aa=\bbbk[t^{\pm 1}]$ be the algebra of Laurent polynomials in variable $t$, and let $c,d$ be two symbols. Set
$$\tilde\fg=\aff(\fg,\sg)=\sum_{i\in\bbbz}(\fg^{\bar i}\otimes t^i)\oplus\bbbk c\oplus\bbbk d,$$
and
$$\tilde\fh:=\fh^\sg\oplus\bbbk c\oplus\bbbk d.$$
Then $\tilde\fg$ becomes a Lie algebra with the bracket:
\begin{eqnarray*}
	[x\otimes t^n+rc+sd,y\otimes t^{n'}+r'c+s'd]&=&[x,y]\otimes t^{n+n'}
	+n\d_{n+n',0}(x,y)c\\
	&&+sn' y\otimes t^{n'}-s'n x\otimes t^n.
	\end{eqnarray*}
The form on $\fg$ extends by linearity to a form on $(\fg\otimes\aa)\oplus\bbbk c\oplus\bbbk d$ as follows:
$$(x\otimes t^n+rc+sd,y\otimes t^{n'}+r'c+s'd)=(x,y)\d_{n+n',0}+rs'+r's.$$
This induces a symmetric invariant bilinear form on $\tilde\fg$ by restriction, which we denote it again by $\fm$. 	 

Henceforth, {\it we assume that}
$$\hbox{there exists $\a\in R^\times$ with  $(\pi(\a),\pi(\a))\not=0$},$$ 
as this condition is always satisfied for the algebras of interest under a suitable choice of 
$\sg$.

\begin{thm}\label{thm:abp1}\cite{ABP02}
	Let $(\fg,\fm,\fh)$ be an extended affine Lie algebra with root system $R$, and $\sg$ an automorphism of $\fg$ satisfying (\ref{eq:axiom}). 
	Then 
	
	(i) $(\tilde\fg,\fm,\tilde\fh)$ is an extended affine Lie algebra,
	
	(ii) if $\fg$ is tame, then so is $\tilde\fg$,
	
	(iii) if $m$ is prime, then the condition $C_{\fg^\sg}(\fh^\sg)=\fh^\sg$ is equivalent to the condition that if
	$\a\in R$ is non-zero then $\pi(\a)\not=0$.
	\end{thm}
	
Considering Theorem \ref{thm:abp1}, we now describe the root system $\tilde R$ and the root spaces of $\tilde\fg$. We have
$$\tilde\fg=\sum_{\tilde\a\in{\tilde \fh}^\star}{\tilde\fg}_{\tilde\a}\hbox{ with }
{\tilde\fg}_{\tilde\a}=\{x\in\tilde\fg\mid [h,x]=\tilde\a(h)x\hbox{ for all }h\in\tilde\fh\}.$$
 Note that we have ${\tilde\fh}^\star=({\fh^\sg})^\star\oplus\bbbk\lam\oplus\bbbk\d$, where
$$\d(\fh^\sg)=\lam(\fh^\sg)=\{0\},\quad \d(c)=\lam(d)=1,\andd \d(d)=\lam(c)=0.$$
It turns out that
$$\tilde{R}\sub\pi(R)+\bbbz\d,$$
and for $\tilde\a=\pi(\a)+i\d\in\tilde R$,
\begin{equation}\label{eq:3}\tilde{\fg}_{\tilde\a}=\left\{\begin{array}{ll}
		\sum_{\{\b\in R\mid\pi(\b)=\pi(\a)\}}\pi_i(\fg_\b)\otimes t^i,
		&\hbox{if }\tilde\a\not=0\vspace{1mm}\\\\
		\tilde\fh,& \hbox{if }\tilde\a=0.\\
		\end{array}\right.
\end{equation}

\begin{lem}\label{lem:prime1}
	Suppose that $m$, the order of $\sg$, is prime. Then for any non-zero integer $k$, 
	$${\tilde\fg}_{k\d}=\pi_k(\fh)\otimes t^k.$$	Furthermore, if $\fg$ is a finite-dimensional simple Lie algebra, then
	${\tilde\fg}_{k\d}\sub\tilde{\fg}_c$.
	\end{lem}
	
	\proof
	Since $k\d\in\tilde{R}\sub\pi(R)+\bbbz\d\sub\la R\ra\oplus\bbbz\d$, we have
	$k\d=\pi(\a)+i\d$ for some $\a\in R$ and $i\in\bbbz$, implying that $\pi(\a)=0$ and $i=k$. But since $m$ is prime, we obtain from Theorem \ref{thm:abp1}(iii) that $\a=0$. Therefore, by 
	(\ref{eq:3}), ${\tilde\fg}_{k\d}=\pi_k(\fg_0)\otimes t^k=\pi_k(\fh)\otimes t^k.$

	Next, let  $\fg$ be a finite-dimensional simple Lie algebra. Then $\fh=\sum_{\a\in{ R}^\times}\bbbk t_\a$ where $R$ is an irreducible finite root system. Therefore, 
	$\tilde{\fg}_{k\d}=\sum_{\a\in R^\times}\bbbk\pi_k(t_\a)\otimes t^k$. Now for each $\a\in R^\times$, we have $\pi(\a)\not=0$ and 
	\begin{eqnarray*}
		\bbbk\pi_{k}(t_\a)\otimes t^k&=&\pi_k[\fg_\a,\fg_{-\a}]\otimes t^k\\
		&=&	\sum_{i}\pi_k [\fg_\a,\pi_{i}(\fg_{-\a})]\otimes t^k\\
(\hbox{by (\ref{eq:4})})&=&	\sum_{i} [\pi_{k-i}(\fg_\a),\pi_i(\fg_{-\a})]\otimes t^k\\
&=&	\sum_{i} [\pi_{k-i}(\fg_\a)\otimes t^{k-i},\pi_i(\fg_{-\a})\otimes t^i]\\
	&\sub&
\sum_{i}	[\tilde{\fg}_{\pi(\a)+(k-i)\d},\tilde{\fg}_{\pi(-\a)+i\d}]\sub{\tilde\fg}_c.
	\end{eqnarray*}	
	\qed

		\begin{lem}\label{lem:pi1}
		Consider distinct roots $\a,\b\in R$ with $\pi(\a)=\pi(\b)$. Then 
		$\pi_k[\fg_\a,\fg_{-\b}]\in C_\fg(\fh^\sg)$ for each $k$.
		Moreover, $\pi[\fg_\a,\fg_\b]=0$.
	\end{lem}
	
	\proof  Let $x\in\fg_\a$ and $y\in\fg_{-\b}$. Then $[x,y]\in\fg_{\a-\b}$. We may assume that $\a-\b$ is a root since otherwise we are done by triviality.
	Now for each $h\in\fh^\sg$, we have
	$$
		[h,[x,y]]=(\a-\b)(h)[x,y]=\pi(\a-\b)(h)[x,y]=0,
	$$
	that is $[x,y]\in C_\fg(\fh^\sg)$. Since $C_\fg(\fh^\sg)=\sum_{\bar i\in\bbbz_m}\fg^{\bar i}\cap C_{\fg}(\fh^\sg)$, we get 
	$$\pi_k[x,y]\in C_\fg(\fh^\sg).
	\qquad\qquad\qquad(\star)
	$$
	If $k=0$, then by $(\star)$ and (\ref{eq:axiom}), we have 
	$$\pi_0[\fg_\a,\fg_\b]=\pi[\fg_\a,\fg_\b]\in\fg^\sg\cap C_\fg(\fh^\sg)=C_{\fg^\sg}(\fh^\sg)\sub\fh^\sg\
	\sub\fg_0.	\qquad\qquad\qquad(\star\star)$$
		On the other hand for each $i$, $\sg^i[x,y]\in\fg_{\sg^i(\a-\b)}$, 
		where $\sg^i({\a-\b})$ is a non-zero root. Thus $\pi[x,y]$ resides in the space spanned by non-zero root spaces.  This together with $(\star\star)$ gives
	$\pi[x,y]=0.$\qed

Let $R_\pi^\times=\{\a\in R^\times\mid (\pi(\a),\pi(\a))\not=0\}$. Since $\sg$ preserves the form, the cyclic group
$\la\sg\ra$ acts on $R^\times_\pi$. We denote by ${\mathcal O}_\sg(\a)$ the orbit of $\a\in R^\times_\pi$. We choose  ${\mathcal O}_\sg$ 
to be a complete set of orbit representatives for this action; therefore $|{\mathcal O}_\sg\cap{\mathcal O}_\sg(\a)|=1$ for each $\a\in R\pi^\times$.
Let  $n_{\sg,k}$ be the number of $\a\in{\mathcal O}_\sg$ for which $\pi_k(\a)\not=0$.

	\begin{rem}\label{rem:1}
	It follows from (\ref{eq:4}) that if $\a$ and $\b$ belong to the same orbit of the action of $\la\sg\ra$ on $R_\pi^\times$, then $\pi(\a)=\pi(\b)$. Conversely, it is known that if $\sg$ is an automorphism arising in the affinization (or twisting) process of a finite-dimensional simple Lie algebra (see \cite{Kac90}) or an affine Kac--Moody Lie algebra (see \cite{Pol94}), then the converse also holds: if $\pi(\a)=\pi(\b)$, then $\a$ and $\b$ lie in the same orbit. In other words, $\pi$ {\it separates} the orbits of the action of $\la\sg\ra$ on $R_\pi^\times$. 
\end{rem}
	
	\begin{pro}\label{pro:2}
	Let $0\not=k\in\bbbz$ and $\pi(\a)+r\d, -\pi(\a)+s\d\in {\tilde R}^\times$, where
	$r+s=k$.  Then
	$$[\tilde{\fg}_{\pi(\a)+r\d},\tilde{\fg}_{-\pi(\a)+s\d}]=\sum_{\{\b \in R_\pi^\times\mid \pi(\b)=\pi(\a)\}}\bbbk\pi_{k}(t_\b)\otimes t^k,\quad(\mod\;\;A)
$$
where $A:=C_\fg(\fh^\sg)\cap\big(\sum_{\mu\in R\setminus\{0\},\pi(\mu)=0}\pi_k(\fg_\mu)\otimes t^k\big)$. 
Moreover, if $\pi$ separates the orbits of $\la\sg\ra$ acting on $R^\times$,
then  
$$\begin{array}{cc}
	[\tilde{\fg}_{\pi(\a)+r\d},\tilde{\fg}_{-\pi(\a)+s\d}]=\bbbk\pi_{k}(t_\a)\otimes t^k &(\mod\; A),\vspace{2mm}\\ \tilde{\fg}_{k\d}\cap\tilde{\fg}_c=\sum_{\a\in{\mathcal O}_\sg}\bbbk t_{{\pi_k(\a)}}\otimes t^k &(\mod\; A),
	\end{array}
	$$
furthermore, if $m$, the order of $\sg$, is prime then $A=\{0\}$ and
	$$\dim (\tilde{\fg}_{k\d}\cap\tilde{\fg}_c)=\dim\sum_{\a\in{\mathcal O}_{\sg,k}}\bbbk\pi_k(\a)\leq \min(\dim\pi_k(\fh_c),\frac{1}{2}n_{\sg,k}).$$ 
	\end{pro}
	
	\proof By (\ref{eq:3}), we have (in the below computations sums run
	over $\b,\gamma\in R_\pi^\times$ with $\pi(\b)=\pi(\gamma)=\pi(\a)$),
	\begin{eqnarray*}
		[\tilde{\fg}_{\pi(\a)+r\d},\tilde{\fg}_{-\pi(\a)+s\d}]
	&=&\sum_{\b}\sum_{\gamma}[\pi_{r}({\fg}_{\b})\otimes t^r,\pi_s({\fg}_{-\gamma})\otimes t^s]\\
(\hbox{by (\ref{eq:4}})&	=&\sum_{\b,\gamma}\pi_{k}[\fg_\b,\pi_s(\fg_{-\gamma})]\otimes t^k\\
(\hbox{by (\ref{eq:pi1}})&=&	\sum_{\b,\gamma}\sum_{i}\pi_{k}[\fg_\b,\omega^{-is}\sg^i(\fg_{-\gamma}]\otimes t^k\\
&=&	\sum_{\b,\gamma}\sum_{i}\pi_{k}[\fg_\b,\fg_{\sg^{i}(-\gamma)}]\otimes t^k\\
(\hbox{by Lemma \ref{lem:pi1}})&=&	\big(\sum_{\b}\pi_{k}[\fg_\b,\fg_{-\b}]+\sum_{\b,\gamma}\sum_{\{i\mid\b-\sg^i(-\gamma)\not=0\}} C_\fg(\fh^\sg)\cap\pi_k(\fg_{\b-\sg^i(-\gamma)})\big)\otimes t^k\\
&=&	\sum_{\b}\bbbk\pi_{k}(t_\b)\otimes t^k\quad(\mod\; A).
\end{eqnarray*}
		We note that the last equality follows from the fact that
		$\pi(\b-\sg^i(-\gamma))=\pi(\a)-\pi(\a)=0$.
	This completes the proof of the first assertion.

Next, assume that $\pi$ separates the orbits of the action. So if $\pi(\b)=\pi(\a)$, we have  $\a,\b$ belong to the same orbit, that is, $\b=\sg^i(\a)$ for some $i$. Then
by (\ref{eq:4}), $\pi_k(\b)=\pi_k(\sg^i(\a))=\omega^{ik}\pi_k(\a)$. 
Thus 
$$\sum_{\{\b\in R_\pi^\times\mid\pi(\b)=\pi(\a)\}}\bbbk\pi_{k}(t_\b)\otimes t^k=
\sum_{\b\in{\mathcal O}_\sg(\a)}\bbbk\pi_{k}(t_\b)\otimes t^k=
\bbbk t_{\pi_k(\a)}\otimes t^k.$$ 
Combining this with the first part of the statement, we conclude that
$$\tfg_{k\d}\cap\tfg_c=\sum_{r+s=k}\sum_{\a\in R_\pi^\times}
	[\tilde{\fg}_{\pi(\a)+r\d},\fg_{-\pi(\a)+s\d}]=
\sum_{\a\in{\mathcal O}_{\sg,k}}\bbbk t_{{\pi_k(\a)}}\otimes t^k\quad(\mod\; A).$$
Finally, if $m$ is prime, then by Theorem \ref{thm:abp1}(iii), for all $\mu \in R \setminus \{0\}$, we have $\pi(\mu) \neq 0$, so $\{i \mid \beta - \sigma^i(-\gamma) \neq 0\} = \emptyset$ and thus $A = \{0\}$. Using $t_{-\alpha} = -t_\alpha$, the last assertion follows.\qed

\begin{exa}\label{exa:1} In this example, we apply Proposition \ref{pro:2} to the finite dimensional simple Lie algebra $\fg$ of type $D_4$ with the root system
	$R=\{\pm(\ep_i\pm\ep_j)\mid 1\leq i,j\leq 4\}$, where $\ep_i$'s are the standard basis of $\bbbr^4$, and $\sg$ is the diagram automorphism of $\fg$ of order $3$. 
	
	Consider the the root base  
	$\{\a_1=\ep_1-\ep_2,\a_2=\ep_2-\ep_3,\a_3=\ep_3-\ep_4,\a_4=\ep_3+\ep_4\}$ for $R$. Then $\sg$ is induced by 
	$\a_1\mapsto\a_3\mapsto\a_4\mapsto\a_1,\andd\a_2\mapsto\a_2.
	$
	Table 1 lists the orbits of the action of $\la \sg\ra$ on positive roots of $D_4$. 
	\begin{table}[H]\label{tab:1}
		\centering
		\caption{}
		\vspace{-10mm}
		$$
		\begin{array}{|c|}
			\hline
			\text{Orbits} \\
			\hline
			{\mathcal O}_1 = \{\epsilon_1 - \epsilon_2, \epsilon_3 - \epsilon_4, \epsilon_3 + \epsilon_4\}  \\
			{\mathcal O}_2 = \{\epsilon_1 - \epsilon_3, \epsilon_2 - \epsilon_4, \epsilon_2 + \epsilon_4\}  \\
			{\mathcal O}_3 = \{\epsilon_1 - \epsilon_4, \epsilon_1 + \epsilon_4, \epsilon_2 + \epsilon_3\} \\
			{\mathcal O}_4 = \{\epsilon_1 + \epsilon_2\} \\
			{\mathcal O}_5 = \{\epsilon_1 + \epsilon_3\}  \\
			{\mathcal O}_6 = \{\epsilon_2 - \epsilon_3\} \\
			\hline
		\end{array}
		$$
	\end{table}
	
	From Table 2 and Proposition \ref{pro:2}, we see that $\pi$ separates the orbits and 
	\begin{table}[H]\label{tab:2}
		\centering
		\caption{}
		\vspace{-10mm}
		$$
		\begin{array}{|c|c|c|c|}
			\hline
			{\mathcal O}_\sg& \pi(\alpha)&\pi_1(\a)&\pi_2(\a) \\
			\hline
			\a_1 = \epsilon_1 - \epsilon_2 & \pi(\alpha_1)&\pi_1(\a_1)&\overline{\pi_1(\a_1)} \\
			\a_1+\a_2 =\epsilon_1 - \epsilon_3 & \pi(\alpha_1) + \alpha_2&\pi_1(\a_1)& \overline{\pi_1(\a_1)}\\
			\a_1+\a_2+\a_3 = \epsilon_1 - \epsilon_4 & 2\pi(\alpha_1) + \alpha_2 &(1+\omega)\pi_1(\a_1)&\overline{(1+\omega)\pi_1(\a_1)}\\
			\a_1+2\a_2+\a_3+\a_4=\epsilon_1 + \epsilon_2 & 3\pi(\alpha_1) + 2\alpha_2&0&0 \\
			\a_1+\a_2+\a_3+\a_4=\epsilon_1 + \epsilon_3 & 3\pi(\a_1)+\alpha_2&0& 0\\
			\a_2 = \ep_2-\ep_3 & \alpha_2 &0&0 \\
			\hline
		\end{array}
		$$
	\end{table}
	$$\dim\tfg_{ (2k+1)\d}=\dim\tfg_{ (2k+2)\d}=1\andd \dim\tfg_{ 3k\d}=2.$$
	In fact 
	$$\begin{array}{l}
		\span_\bbbk\{ \pi_1(\a)\mid \a\in{\mathcal O}_{\sg,1}\}=\bbbk\pi_1(\a_1),\\
		\span_\bbbk\{ \pi_2(\a)\mid \a\in{\mathcal O}_{\sg,2}\}=\bbbk\overline{\pi_1(\a_1)},\\
		\span_\bbbk\{ \pi_0(\a)\mid \a\in{\mathcal O}_{\sg,0}\}=\bbbk\pi(\a_1)\oplus\bbbk\a_2.
	\end{array}
	$$
	(Note that an expression of the form $\overline{\sum_i{k_i\a_i}}$, $k_i\in\bbbk$, is defined to mean $\sum_i\bar{k}_i\a_i$.
\end{exa}

\begin{cor}\label{cor:8-1}
	Suppose that $m$ is prime and $k$ is a non-zero 
integer. Then
$${\tilde\fg}_{k\d}=\pi_k(\fh)\otimes t^k.$$	Furthermore, if $\fg$ is a finite-dimensional simple Lie algebra, then
${\tilde\fg}_{k\d}\sub\tilde{\fg}_c$.
\end{cor}

\proof	Since $k\d\in\tilde{R}\sub\pi(R)+\bbbz\d\sub\la R\ra\oplus\bbbz\d$, we have
$k\d=\pi(\a)+i\d$ for some $\a\in R$ and $i\in\bbbz$, implying that $\pi(\a)=0$ and $i=k$. But since $m$ is prime, we obtain from Theorem \ref{thm:abp1}(iii) that $\a=0$. Therefore, by 
(\ref{eq:3}), ${\tilde\fg}_{k\d}=\pi_k(\fg_0)\otimes t^k=\pi_k(\fh)\otimes t^k.$

Next, let  $\fg$ be a finite-dimensional simple Lie algebra. Then $\fh=\sum_{\a\in{ R}^\times}\bbbk t_\a$ where $R$ is an irreducible finite root system.  Now since $m$ is prime, we have $A=0$, and for  $\a\in R^\times$, we have have $\pi(\a)+i\d\in{\tilde R}^\times$. Thus using Proposition \ref{pro:2}, we obtain  
\begin{eqnarray*}
\tilde{\fg}_{k\d}=\sum_{\a\in R^\times}\bbbk\pi_k(t_\a)\otimes t^k
	&\sub&
\sum_{\a\in R^\times}	\sum_{i}	[\tilde{\fg}_{\pi(\a)+(k-i)\d},\tilde{\fg}_{\pi(-\a)+i\d}]\sub{\tilde\fg}_c.
\end{eqnarray*}	
\qed

\begin{thm}\label{thm:6}
Let $(E,\fm,\hh)$ be an extended affine Lie algebra of of nullity $1$ with root system $R.$ Then $(\hat E,\fm,\hat\hh)$ is a an affine Lie algebra.
\end{thm}

\proof  By Proposition \ref{pro:siz2} and Lemma \ref{lem:5}, we only need to show the form $\fm$ on
$\hat{E}^{iso}$ is non-degenerate.

Let $E$ be of type $X$. Since $E$ has nullity $1$, the root system $R$ of $E$  is an affine root system, that is, $R=\bbbz\d\cup(\rd+\bbbz\d)$ where $\rd$ is a finite root system of type $X$ and $\d$ is a non-zero isotropic root.
We must show that for each non-zero isotropic root $\lam=k\d$,  the form $\fm$ on $(E_\lam\cap E_c)\oplus (E_{-\lam}\cap E_c)$ is non-degenerate. 

By \cite[Proposition 1.28]{AG01} (or \cite[Theorem 7.3]{Yos06}), 
the algebra $E_c/Z(E_c)$ is a centerless Lie torus of type $X$ of nullity $1$. 
Thus we may identify it with the core modulo center of an affine Lie algebra of type $X$. 
Now, one knows that any affine Lie algebra is obtained by the affinization (twisting) process from a finite dimensional simple Lie algebra $(\dot\fg,\fm,\dot\fh)$ whose root system can be identified with $\rd$, and a (diagram) automorphism $\sg$ of prime order.
Therefore, modulo $Z(E_c)$, we may identify $E_{k\d}\cap E_c$ with 
$\tilde{\dot\fg}_{k\d}\cap\tilde{\dot\fg}_c$ which is equal to $\tilde{\dot\fg}_{k\d}$, by Corollary \ref{cor:8-1}. 
Then each non-zero element in $E_{k\d}\cap E_c$ has the  form 
$x+z$ where $0\not=x\in\tilde{\dot\fg}_{k\d}$, and
$z$ in 
$Z(E_c)$. By  By (\ref{eq:khat1}),  there exists $y\in \tilde{\dot\fg}_{-k\d}\equiv E_{-k\d}\cap E_c$, $\mod Z(E_c)$,  such that
$(x,y)\not=0$. 
Since $(Z(E_c),E_c)=\{0\}$ (see (\ref{eq:9}), we have $(x+z,y)\not=0$. 

\qed

\begin{cor}\label{cor:10} (Affine Localization)
	Let $(E,\fm, R)$ be an extended affine Lie algebra with root system $R$. Let $R'$ be an affine root subsystem of $R$. Then $E_{R',\hh_{R'}}$ 
		is an extended affine Lie algebra of nullity $1$, and $\big(\hat{E}_{R',\hh_{R'}},\fm,\hat{\hh}_{R'}\big)$ is an affine Lie subalgebra of $E$.
\end{cor}

\proof  
By \pref{parag:n1}, $R'$ is a closed subsystem of $R$. Therefore, the first assertion is an immediate consequence of Proposition \ref{pro:khan2}.	Then Theorem \ref{thm:6} implies that the subalgebra $\big(\hat{E}_{R',\hh_{R'}},\fm,\hat{\hh}_{R'}\big)$ of $E$ 
is an affine Lie algebra.\qed

Recall from (\pref{pt1}) that $R\sub\rd +\Lam$. Let $X$ be the type of $R$ (or equivalently the type of $\rd$). Let $0\not=\d\in R^0$. As in Example \ref{exa:-1}, we set 
$$R_{\dot R,\d}:=R^\times_{\dot R,\d}\cup R^0_{\dot R,\d},$$
where
$$
R^\times_{\dot R, \d}:=(\rd+\bbbz\d)\cap R^\times\hbox{ and } R^0_{\dot R,\d}=(R^\times_{\dot R,\d}-R^\times_{\dot R,\d})\cap R^0.$$
	Also,  in accordance with the notation of Section \ref{sec:local}, we set
	$$\hat{\hh}_{\dot R,\d}:=\hh_{R_{\dot R,\d}}\andd\hat{E}_{\dot R,\d}:=\hat{E}_{R_{\dot R,\d},\hat{\hh}_{\dot R,\d}}.
	$$  
	Then from Corollary \ref{cor:10}, we have the following.

\begin{cor}\label{cor:11}
	Let $(E,\fm, R)$ be an extended affine Lie algebra with root system $R$. Let $0\not=\d\in R^0$. Then  $(\hat{E}_{\dot R,\d},\fm,\hat{\hh}_{\dot R,\d})$  is an affine Lie subalgebra of $E$ of the same type of $E$.
	\end{cor}

	\begin{cor}\label{cor:12} 
		Let $(E,\fm,\hh)$ be an extended affine Lie algebra. Let $\a\in R^\times$, $\sg\in R^0\setminus\{0\}$ with $\a+\sg\in R^\times$. Then the root spaces  $E_{\pm\a}$, $E_{\pm(\a\pm\sg)}$
		 can be embedded in 
		 an affine Lie subalgebra of $E$ of 
		 type $A_1$. 
	\end{cor}
	
	\proof 
	Set $T=\{\a,\a+\sg\}$. According to \pref{parag:n1} and Example \ref{exa:-1}, $\tilde{R}_T={R}_{\a,\sg}$ is a closed  affine subsystem of $R$ of type $A_1$.
	Then $\hat{E}_{R_{\a,\sg}},\hat{\hh}_{R_{\a,\sg}}$ 
		is an affine subalgebra of $E$ containing the root spaces
		$E_{\pm\a}$, $E_{\pm(\a\pm\sg)}$.
	\qed 
	
To see an application of affine localization in the study of Chevalley bases for EALAs,  we recall the concept of a Chevalley system for the EALA $E$. Let $\tau$ be an automorphism of $E$ of order $2$  (an involution ). A family
$\{x_\a\}_{\a\in R^\times}$ is called a {\it Chevalley system} 
	for $E$ with respect to $\tau$ if for each $\a\in R^\times$,
	$x_\a\in E_\a$, $(x_\a,[x_\a,x_{-\a}],x_{-\a})$ is an 
	$\frak{sl}_2$-triple, and
	$\tau(x_\a)=-x_{-\a}$.
	
	\begin{cor}\label{coe:az1}
		Let $\{x_\a\}_{\a\in R^\times}$ be a Chevalley system for $E$ with respect to an involution.
		Let $\a,\b\in R^\times$ and $\d\in R^0\setminus\{0\}$ such that $\a+\d,\b+\d\in R^\times$ and $\a\pm\b$ is not isotropic. Then
		$[[x_{\a+\d},x_{-\a}],x_\b]\in2\bbbz x_{\b+\d}$.
		\end{cor}
		
		\proof
	If both $\b-\a$ and $\b+\a+\d$ are not roots, then the Jacobi identity implies
			$[[x_{\a+\d},x_{-\a}],x_\b]=0$, and we are done. Therefore, we may assume that at least one of $\b-\a$ or $\b+\a+\d$ is a root. Consequently, depending whether  $(\a,\b)\not=0$ or $(\a,\b)=0$, at least one of the following is a connected subset of $R^\times$:
			$$\{\a,\b,\a+\d\},\quad\{\a,\a+\b+\d,\b\},\quad \{\b-\a,\a,\a+\d\}.$$
			We denote one such connected subset by $T$. Then the triple
			$(\tilde{R}_T,\fm_T,{\tilde\v}_T)$ is an EARS of nullity $1$ (see  \pref{parag:n1}), and hence an affine root system.  
	
		By Corollary \ref{cor:11}, the triple $(\hat{E}_{\tilde{R}_T},\fm,\hat{\hh}_{\tilde{R}_T})$ 
			is an affine Lie subalgebra of $E$ that contains the root spaces $E_{\pm\a}$, $E_{\pm(\a+\d)}$, $E_{\pm\b}$ and $E_{\pm(\b+\d)}$. Thus the Lie bracket in the statement can be evaluated within this local affine Lie subalgebra. The result then follows from the known realization of affine Lie algebras; see \cite{Mit85}.
		\qed	
	 
	We conclude this section with a result concerning ``finite localization''.
	
	\begin{pro}\label{pro:11} (Finite Localization)
		Let $R'$ be a finite subsystem of $R$. Let $\hh'$ be a cover for $R'$ in $\hh$.
		Then the Lie cover $E':=E_{R',\hh'}$ is finite dimensional reductive Lie algebra.
		In particular $E_{R',\hh_{R'}}$ is a finite dimensional simple Lie algebra with root system $R'$.	
		\end{pro}
	
	\proof By \pref{parag:n1}, $R'$ is a closed subsystem of $R$. Since a finite extended affine root system is an irreducible finite root system, from Proposition \ref{pro:khan2} we see that  $E_{R',\hh'}$ is an extended affine Lie algebra of nullity $0$. Thus it is a reductive Lie algebra. Now since $\dim\hh_{R'}=\dim\sum_{\a\in R'}\bbbk t_\a=\rank\, R'$, we conclude that  $E_{R',\hh_{R'}}$ is a finite dimensional simple Lie algebra.
	\qed

\section{Examples}\setcounter{equation}{0}\label{examples}

In Section \ref{sec:affinization}, we proved a systematic procedure of constructing affine Lie subalgebras of an extended affine Lie algebra. We now apply the results of Section \ref{sec:affinization} to some specific extended affine Lie algebras and explore the corresponding affine Lie subalgebras. In particular, we want to see the structure of the non-zero isotropic root spaces.
     
	 \begin{exa}\label{exa:2}
		(Toroidal extended affine Lie algebras). Let $(\dot\fg,\fm,\dot\hh)$ 
		be a finite dimensional simple Lie algebra with root system $\rd$. 
		Let $\aa=\bbbk[t_1^{\pm 1},\ldots, t^{\pm 1}_\nu]$ be the ring of Laurent polynomials in $\nu$ variables. 
		Let $\Lam=\bbbz^\nu$. For $\lam=(n_1,\ldots,n_\nu)\in\Lam$, we set $t^\lam:=t_1^{n_1}\cdots t_\nu^{n_\nu}$. Then $\aa=\sum_{\lam\in\Lam}\aa^\lam$, where $\aa^\lam=\bbbk t^\lam.$ We have the root space decomposition  $\dot\g=\dot\fh+\sum_{\dot\a\in\rd^\times}{{\dot\fg}_{\dot\a}}$,  with respect to a Cartan subalgebra $\dot\fh$. Set
		$R=\dot R+\Lam$ and $E=\sum_{\a\in R}E_\a$, where $E_0=\hh:=\dot\fh\oplus D\oplus C$ and for $\a=\dot\a+\lam$, $\dot\a\in \rd^\times$, $\lam\in \Lam$,
		$E_\a={\dot\fg}_{\dot\a}\otimes\aa^\lam$.	Here  $D=\sum_{i=1}^\nu\bbbk d_i\sub\hbox{Der}(\aa)$ and $C=D^\star$ where 
		$d_i=t_i\frac{\partial}{\partial t_i}$. Each $d_i$ can be extended to a derivation of $\dot\fg\otimes\aa$ by $d_i(a\otimes t^\lam)=\lam_i a\otimes t^\lam$. Then $C=\sum_{i=1}^\nu\bbbk c_i$ with $c_i(d_j)=\d_{ij}$.
		The bracket on $\dot\fg$ extends to a bracket on $E$ by 
		$$\begin{array}{c}
		[C,E]=\{0\},\\
		\;[d,a\otimes t^\lam]=a\otimes d(t^\lam),\\	
	\;[a\otimes t^\lam,a'\otimes t^{\lam'}]=[a,a']\otimes t^{\lam+\lam'}+\d_{\lam,-\lam'}(a,a')\sum_{i=1}^\nu\lam_ic_i,
		\end{array}
		$$
		and 
		the form on $\dot\fg$ extends to $E$ by
		$$(a\otimes t^\lam+c+d,a'\otimes t^{\lam'}+c'+d')=\d_{\lam,-\lam'}(a,a')+c(d')+c'(d).$$ 
		Here $a,a'\in\dot\fg$, $c,c'\in C$, $d,d'\in D$, $\lam,\lam'\in\Lam$. Then the triple
		$(E,\fm,\hh)$ is an extended affine Lie algebra with root system $R$. The type of $E$ is the type of $\dot\fg$. In the notation of Section \ref{sec:tame6}, we have $E=E(\fg,D,\kappa)$, where 
		$$\fg=\dot\fg\otimes\aa,\; D=\mathrm{SCDer}(\fg)^0=\{\partial_\theta\mid\theta\in\Hom_\bbbz(\Lam,\bbbk)\},\hbox{ and }\kappa=0.$$ Note that, we have $E_c=\fg\oplus C$, and $Z(E_c)=C=E_c^\perp.$ This shows that $E$ is tame.
		
		Next for $0\not=\d\in\Lam$, consider the affine root system 
		$R_{\dot R,\d}$, see Example \ref{exa:-1}.
		By Corollary \ref{cor:11},  the subalgebra $(\hat{E}_{\dot R,\d},\fm,\hat{\hh}_{\dot R,\d})$ is an affine Lie subalgebra of type $\dot R$. 
		
		Note that if $k\not=0$, the root space $\hat{E}_{\dot R,\d}$, corresponding to the root $k\d$, is spanned by elements of the form $\sum_{r+s=k}[E_{\a+r\d},E_{-\a+s\d}]$, $\a\in R^\times_{\dot R,\d}$.
			Now for $\a\in R^\times$, we have $\a=\dot\a+u\d$, where $\dot\a\in\rd^\times$, $u\in\bbbz$. Then
		\begin{eqnarray*}
			\sum_{r+s=k}[E_{\a+r\d},E_{-\a+s\d}]&=&\sum_{s\in\bbbz}[E_{\dot\a+(u+s+k)\d},E_{-\dot\a-(u+s)\d}]\\
			&=&\sum_{s\in\bbbz}[{\dot\fg}_{\dot\a}\otimes t^{(u+s+k)\d},{\dot\fg}_{-\dot\a}\otimes t^{(-u-s)\d}]\\
			&=&\bbbk t_{\dot\a}\otimes t^{k\d}.
		\end{eqnarray*}
	Thus $(\hat{E}_{\dot R,\d})_{k\d}=\sum_{\dot\a\in\rd}\bbbk t_{\dot\a}\otimes t^{k\d}$.
		\end{exa}
	
 \begin{exa} (Type $B_\ell$)
  	We begin by a general setup, see \cite[III.\S3]{AABGP97}. Let $\nu \geq 1$ and
 	$\aa$ be the algebra of Laurent polynomials in variables $t_1,\ldots, t_\nu$, namely
  	$\aa = \sum_{\sigma\in\Bbb{Z}^\nu}\bbbk t^\sigma$.
 	
 	Let
 	$\Lambda := \zn$.  
 	Suppose that $m\ge 1$ and $\tau_1,\dots,\tau_m$ represent distinct cosets
 	of  $2\zn$ in $\zn$ with
 	$
 	\tau_1 = 0$. Then $S =
 	\cup_{i=1}^m(\tau_i+2\zn)$ is a semilattice in $\Bbb
 	R^\nu$ such that
 	$2\zn \sub S \sub  \zn.$
 	
 	Now let
 	$n=\ell+m$ where $\ell \ge 2$. Put
 	$$F = \left[\begin{array}{lll}
 		t^{\tau_1} &\cdots & 0 \\
 		\vdots & \ddots & \vdots\\
 		0 & \cdots & t^{\tau_m} \end{array}\right]\quad \text{ and
 	}\quad
 	G =\left[ \begin{array}{lll}
 		0 & I_\ell & 0 \\
 		I_\ell & 0 & 0 \\
 		0 & 0 & F \end{array}\right].
 	$$
  	Consider the subalgebra
 	$$\fg :=
 	\{ X\in M_n(\aa): G^{-1}X^tG=- X\hbox{ and }\tr (X)\equiv 0\hbox{ mod }[\aa,\aa]\}.$$
 	of $\gl_n(\aa)$. One endows $\fg$ with a symmetric invariant non-degenerate form 
 	by linear extension of 
 	$$
 	(t^{\mu} e_{i,j}, t^{\nu} e_{r,s}) = \delta_{\mu, -\nu} \delta_{j,r} \delta_{i,s},
 	$$
 	where $\{e_{i,j}\mid 1\leq i,j\leq n\}$ is the set of matrix units.

 	Next, consider the $\ell$-dimensional abelian subalgebra
 	$\dot{\fh}$ of $\fg$ spanned by matrices $e_{ii}-e_{\ell+i,\ell+i}$, $1\leq i\leq \ell$.
  Define
 	$\ep_i\in \dot{\fh}^*$ by
 	$\ep_i
 	(e_{jj}-e_{\ell+j,\ell+j})=\d_{ij},$
 	for $i =1, \dots, \ell.$
 	Then,
 	$$
 	\fg= \sum_{\da \in \dot R} \fg_{\dot\a}
 	= \fg_0 +
 	\sum_{1\le i \ne j\le\ell}\fg_{\ep_i\pm\ep_j} 
 	+\sum_{1\le i\le\ell}
 	\fg_{\pm\ep_i} 
 	$$
 	where 
	$$\fg_{\dot\a} = \{ x\in \fg \mid [h, x] = \da(h)x, \text{ for all }
 	h\in \dot{\fh} \}\andd\rd =\{ \da \in \dot{\fh}^* \mid \fg_{\dot\a} \neq
 	\{0\} \}.$$
 	Moreover,
 	\begin{equation}\label{eq:rok4}
 		\begin{array}{c}
 			\fg_{\ep_i-\ep_j} = \{a(e_{ij}-
 			e_{\ell+j,\ell+i}) \mid a\in\aa\},\vspace{2mm}\\
 			\fg_{\ep_i+\ep_j} = \{a(e_{i,\ell+j}-
 			e_{j,\ell+i})\mid a\in\aa\},\vspace{2mm}\\
 			\fg_{-\ep_i-\ep_j} = \{a(e_{\ell+i,j}-
 			e_{\ell+j,i}) \mid a\in\aa\},\vspace{2mm}\\
 			\fg_{\ep_i} = \{\sum_{j=1}^m
 			a_j (e_{2\ell+j,\ell+i} - 
 			x^{\tau_j}e_{i,2\ell+j}),     \mid a_1,\dots,a_m \in\aa\},\vspace{2mm}\\
 			\fg_{-\ep_i} = \{\sum_{j=1}^m
 			a_j (e_{2\ell+j,i} - x^{\tau_j}
 			e_{\ell+i,2\ell+j})
 			\mid a_1,\dots,a_m \in\aa\}.
 			%
 		\end{array}
 	\end{equation}
  	Hence, $\rd=\{\ep_i,\ep_i\pm\ep_j\mid 1\leq i,j\leq\ell\}$ is an irreducible finite root system of type $X$ in $\vd =
 	\sum_{i=1}^\ell \Bbb R \ep_i$, where $X=B_\ell$.
 
  	Next, $\fg$ becomes a $\Lam$-graded Lie algebra by setting
 	$\deg(x^\sigma e_{pq}) = 2\sigma + \lambda_p-\lambda_q,$
 	where
 	$\lambda_1=0,\dots,\lambda_{2\ell}=0,\lambda_{2\ell+1} =
 	\tau_1,\dots,\lambda_n = \tau_m.$
 	Finally, if 
	\begin{equation}\label{eq:temp-1}
		D=D^0=\dd=\{\partial_\theta\mid\theta\in\Hom_\bbbz(\Lam,\bbbk)\},
		\end{equation}
		 then, considering  (\ref{rok1}), the triple $E=E(\fg,D,0)$ is a tame extended affine Lie algebra with root system  $R=(S+S)\cup (\rds +S)\cup (\rdl+L)$, where 
	$$\rds=\{\pm\ep_i\mid 1\leq i\leq \ell\},\quad\rdl=\{\ep_i\pm\ep_j\mid 1\leq i\not=j\leq\ell\},\andd
 	L = 2\zn.$$
 	
	We now fix $0\not =k\in\bbbz$. Let $\a\in R^\times$ and $0\not=\sg\in R^0$. We want to compute
	$$[E_{\a+n\sg},E_{-\a+n'\sg}]\qquad(\star)$$ such that $n+n'=k$, since $E_c\cap E_{k\sg}$ is spanned by such elements.
	We split the computations into cases where $\a$ is short or long.

 	(I) $\a$ is short. Then $\a=\dot\a+\tau_t+2\lam$ for some $1\leq t\leq m$ and $\lam\in 2\Lam$.
	By symmetry between $\pm\ep_i$, we may assume that $\dot\a=\ep_i$ for some $i$. 
 	Now $\a+n\sg=\ep_i+\tau_t+2\lam+n\sg.$
 	Then
 	$$E_{\a+n\sg}=E_{\dot\a+\tau_t+n\sg+2\lam}=\fg_{\dot\a}\cap\fg^{\tau_t+n\sg+2\lam},$$
 	and
 	$$E_{-\a+n'\sg}=E_{-\dot\a-\tau_t+n'\sg-2\lam}=\fg_{-\dot\a}\cap\fg^{-\tau_t+n'\sg-2\lam}.$$
 	By (\ref{eq:rok4}), the root spaces $E_{\a+n\sg}$ and $E_{-\a+n'\sg}$ are spanned by elements of the form 
 	$$\begin{array}{c}
 		X(\gamma,i,r):=t^\gamma e_{2\ell+r,\ell+i}-t^{\gamma+\tau_r} e_{i,2\ell+r},\vspace{2mm}\\
 		\bar X(\gamma,i,r):=t^{\gamma} e_{2\ell+r,i}-t^{\gamma+\tau_r} e_{\ell+i,2\ell+r},
 	\end{array}
 	$$
 where	$1\leq r\leq m$, and $\gamma$ is appropriately chosen. 

 Now
 	\begin{eqnarray*}
 		&&[X(\gamma,i,r),\bar{X}(\gamma',i,s)]\\
 		&&=t^{\gamma+\gamma'}(-t^{\tau_s}e_{2\ell+r,2\ell+s}-\d_{r,s}t^{\tau_r}e_{i,i}+t^{\tau_r}e_{2\ell+s,2\ell+r}+\d_{r,s}t^{\tau_s}e_{\ell+i,\ell+i}),
 	\end{eqnarray*}
 	with
 	$$k\sg=\deg ([X(\gamma,i,r),\bar{X}(\gamma',i,s)])=\left\{\begin{array}{ll}
 		2\gamma+2\gamma'+2\tau_s&\hbox{if }r=s,\\
 		2\gamma+2\gamma'+\tau_s+\tau_r&\hbox{if }r\not=s.\\
 		\end{array}\right.
 		$$
 		Therefore, if $k\sg\in2\Lam$, then $r=s$. So, we conclude that 
 		\begin{equation}\label{eq:temp1}
 	(\star)={\bbbk}	t^{\mu+\tau_s}(e_{i,i}-e_{\ell+i,\ell+i}),\qquad(k\sg=2\mu+2\tau_s).
 			\end{equation}
 If $k\sg\not\in 2\Lam$, then
$r\not=s$. Here, we note that
\begin{equation}\label{ten1}
	\begin{array}{c}
		2\gamma+\tau_r=\deg(X(\gamma,i,r))=\tau_t+n\sg+2\lam,\vspace{2mm}\\
		2\gamma'+\tau_s=\deg(\bar{X}(\gamma',i,s)=-\tau_t+n'\sg-2\lam,
	\end{array}
\end{equation}
therefore, as far as the sum $n+n'$ is fixed, the indices $r,s$ are uniquely determined by $\tau_t$ and $\sg$; since $\tau_i$'s represent distinct cosets of $2\Lam$ in $\Lam$. In fact, since $s\not=r$, from $k\sg=2\mu+\tau_s+\tau_r$, we have $k=n+n'$ is odd and so we may assume that $n$ is odd and $n'$ is even. This forces $s=t$ and $\tau_r=\tau_t+n\sg$ mod $2\Lam$. Thus 
	\begin{equation}\label{eq:temp2}
	(\star)={\bbbk}	t^{\mu}(t^{\tau_s}e_{2\ell+r,2\ell+s}-t^{\tau_r}e_{2\ell+s,2\ell+r}),\qquad (k\sg=2\mu+\tau_t+\tau_r).
	\end{equation}

 		(II) $\a$ is long. Then $\a=\dot\a+\lam$, where $\dot\a\in\rdl$, $\lam\in 2\Lam$. We may assume $\da=\ep_i+\ep_j$, since the computations for the case $\dot\a=\ep_i-\ep_j$ is similar.
 	By (\ref{eq:rok4}), the root spaces $E_{\a+n\sg}$ and $E_{-\a+n'\sg}$ are spanned by elements of the form  		
	$$X(\gamma,i,j)=t^\gamma (e_{i,\ell+j}-e_{j,\ell+i})\andd 
	\bar{X}(\gamma,i,j)=t^\gamma (e_{\ell+i,j}-e_{\ell+j,i}).$$
 		Then
 		$\deg X(\gamma,i,j)=2\gamma=\deg\bar{X}(\gamma,i,j),$
 		and
 		\begin{eqnarray*}
		[X(\gamma,i,j),
	\bar{X}(\gamma',i,j)]&=&
 		[t^\gamma(e_{i,\ell+j}-e_{j,\ell+i}),t^{\gamma'}(e_{\ell+i,j}-e_{\ell+j,i})]\\
 			&=&
 			t^{\gamma+\gamma'}(-e_{ii}-e_{j,j}+e_{\ell+i,\ell+i}+e_{\ell+j,\ell+j}).
 		\end{eqnarray*}
 	Thus 
 	\begin{equation}\label{eq:temp3}
 		(\star)=\bbbk t^{\mu}(-e_{ii}-e_{j,j}+e_{\ell+i,\ell+i}+e_{\ell+j,\ell+j}),\qquad(k\sg=2\mu).
 	\end{equation}

Considering (I) and (II), we conclude from (\ref{eq:temp1})-(\ref{eq:temp3}) that
$$\dim (E_{\dot R,\sg})_{k\sg}=\left\{\begin{array}{ll}
	\ell&\hbox{if }k\sg\in L,\\
	1& \hbox{if }k\sg\not\in L,
	\end{array}\right.
	$$
 	a well-known result concerning the dimension of isotropic root spaces of an affine Lie algebra with the affine label	$B_\ell^{(2)}$ or $D_{\ell+1}^{(2)}$ in the notation of \cite{MP95} or \cite{Kac90}, respectively.
\end{exa}

\begin{exa}\label{exa:A1} (Type $A_1$)
	Let $\Lam=\bbbz^\nu$ and $S$ be a semilattice in $\Lam=\bbbz^\nu$.
That is, $\la S\ra=\Lam$ and $S=\cup_{i=0}^mS_i$ where
$S_i=\tau_i+2\Lam$ and $\tau_0,\tau_1,\ldots,\tau_m$ represent distinct cosets of $2\Lam$ in $\Lam$ with $\tau_0=0$. 
For $\sg\in S$ consider the symbol $x^\sg$ and set $\jj=\jj_S=\sum_{\sg\in S}\bbbk x^\sg$. By convention, we write $\jj=\sum_{\sg\in\Lam}\bbbk x^\sg$, where we interpret $x^\sg=0$ if $\sg\not\in S$. Then $\jj$ is a Jordan algebra with the multiplication
\begin{equation}\label{eq:rokn2}
	x^\sg\cdot x^\tau=\left\{\begin{array}{ll}
	x^{\sg+\tau}&\hbox{if }\sg,\tau\in S_0\cup S_i,\; 0\leq i\leq m,\\
	0&\hbox{otherwise}.
\end{array}\right.
\end{equation}
In fact by setting $\jj^\sg=\bbbk x^\sg$, $\sg\in\lam$, we see that $\jj$ is a $\Lam$-Jordan torus.

Next, let $L_x$ denote the operator on $\jj$ defined by $L_xy=xy$, for $y\in\jj$, and set 
$$\begin{array}{c}
	\mathrm{Inder}(\jj):=\{L_x\mid x\in \jj\},\\\mathrm{Instrl}(\jj):=\sum_{x\in \jj}L_x+\{\sum_i[L_{x_i},L_{y_i}]\mid x_i,y_i\in\jj\},\\ 
\mathrm{TKK}(\jj):=\jj\oplus\mathrm{Instrl}(\jj)\oplus\bar{\jj},
\end{array}
$$
where $\bar\jj$ is a copy of $\jj$.
Then $\mathrm{Instrl}(\jj)$ and $\mathrm{TKK}(\jj)$ are Lie algebras with the brackets given by
$$\begin{array}{c}
	[L_x+C,L_y+D]=[L_x,L_y]+L_{Cy}-L_{Dx}+[C,D],\\
	\begin{array}{c}
		[x_1+\bar{y}_1+E_1,x_2+\bar{y}_2+E_2]=-E_2x_1+E_1x_2-\overline{\bar E_2 y}+\overline{\bar E_1y_2}\\
		\qquad\qquad\qquad\qquad\qquad\qquad\qquad+x_1\triangle y_2-x_2\triangle y_1+[E_1,E_2],
	\end{array}
	\end{array}
	$$
	for $x,y\in\jj,C,D\in\mathrm{Inder}(\jj)$,
$x_i\in\jj,\bar y_i\in\bar{\jj}$, and $E_i\in\mathrm{Instrl}(\jj)$, where $x\triangle y=L_{xy}+[L_x,L_y]$.
Here $\bar{}:\mathrm{Instrl}(\jj)\rightarrow\mathrm{Instrl}(\jj)$ is an involution define by 
$\overline{L_x+D}=-L_x+D$. The Lie algebra $\mathrm{TKK}(\jj)$ is called TKK {\it Lie algebra} of $\jj$.
Now set
$$\fg:=\mathrm{TKK}(\jj)$$
Then $\fg$ is $\Lam$-graded with
	$\fg^\lam=	\mathrm{TKK}(\jj)^\lam=\jj^\lam\oplus\mathrm{Instrl}(\jj)^\lam\oplus\overline{\jj^\lam},$
for $\lam\in\Lam$, where $\mathrm{Instrl}(\jj)^\lam=L_{\jj^\lam}\oplus\sum_{\mu+\nu=\lam}[L_{\jj^\mu},L_{\jj^\nu}]$.

The Lie algebra $\fg$ can be equipped with a symmetric invariant non-degenerate form as follows.  
Consider the linear map $\ep:\jj\rightarrow\bbbk$ induced by $\ep(1)=1$ and $\ep(x^\lam)=0$ if $\lam\not=0$. Then
$(x,y):=\ep(xy)$ defines a symmetric invariant non-degenerate form on $\jj$. This then defines a form on $\mathrm{Instrl}(\jj)$ by 
$(D,[L_x,L_y])=(Dx,y)$ for $D\in\mathrm{Instrl(\jj)}$, $x,y\in\jj$. Finally, this form extends to a symmetric invariant non-degenerate form on $\fg$ by
$$(x_1+L_{y_1}+D_1+\bar{z}_1, x_2+L_{y_2}+D_2+\bar{z}_2):=(x_1,z_2)+(x_2,z_1)+(y_1,y_2)+(D_1,D_2),$$
for $x_i,y_i,z_i\in\jj$, $D_i\in\mathrm{Instrl}(\jj)$. It turns out that $\fg$ is a centerless $\Lam$-Lie torus. Now we consider the extended affine Lie algebra
$E=E(\fg,D,0)$ constructed in Section \ref{sec:tame6}, where $D$ is given by (\ref{eq:temp-1}), and 
$\hh=\bbbk L_1\oplus C\oplus D$ is the Cartan subalgebra of $E$. In fact
$(E,\fm,\hh)$ is an extended affine Lie algebra of type $A_1$ with root system $R\sub\Lam\cup (\pm\dot\a+\Lam)$, where $\dot\a (L_1)=1$ and $\dot\a (C\oplus D)=\{0\}$. Then for $\a=\dot\a+\sg\in R$, we have
$E_{\dot\a+\sg}=\bbbk x^\sg$, $E_{-\dot\a+\sg}=\bbbk\bar{x}^\sg$, $E_0=\bbbk L_1\oplus C\oplus D$ and 
	$E_\sg=\bbbk L_{x^\sg}+\sum_{\tau\in\Lam}\bbbk [L_{x^{\sg+\tau}},L_{x^{-\tau}}]$ if $\sg\not=0$.

Recall that $S=\cup_{i=0}^m S_i$. We consider more cosets $S_i=\tau_i+2\Lam$, $i=m+1,\ldots, 2^\nu-1$ of $2\Lam$ in $\Lam$ such that $\Lam=\cup_{i=0}^{2^\nu-1}S_i$. We set $\ss:=\{S_i\mid 0\leq i\leq 2^\nu-1\}$ and define a symmetric function $\Gamma:\ss\times\ss\rightarrow\{0,1\}$ by
$$\Gamma(S_i,S_j)=\left\{\begin{array}{ll}
	1&\hbox{if $0\leq i,j\leq m$, $i=j$ or $ij=0,$}\\
	0&\hbox{otherwise.}
	\end{array}\right.
	$$ 
Then, for $0\leq i\not=j\leq m$, we have
$$\Gamma(S_0,S_0)=\Gamma(S_0,S_i)=\Gamma(S_i,S_i)=1\andd \Gamma(S_i,S_j)=0,$$
and 
$$
\Gamma (S_i,S_j)=0,\qquad (S_i+S_j\not\in\ss,\; 0\leq i,j\leq m),$$
where the latter equality holds since if $S_j+S_j\not\in\ss$, then $i,j$ can not be equal and none can be zero.
By convention we set
$$\Gamma(\sg,\tau):=\Gamma(S_i,S_j)=\Gamma(\sg,S_j)\qquad(\sg\in S_i,\tau\in S_j).$$ 
Then
$$ 
x^\sg\cdot x^\tau=\Gamma(\sg,\tau)x^{\sg+\tau}=\Gamma(S_i,S_j)x^{\sg+\tau}\qquad(\sg\in S_i,\tau\in S_j,\;0\leq i,j\leq 2^\nu-1).$$
Here, we list some useful facts from \cite[\S 5]{Az25} which are of our use in the sequel. For $i,j,k$, we have
\begin{equation}\label{eq:rokn3}
	\begin{array}{c}
	\Gamma(S_0, S_0+S_i)=\Gamma(S_i,S_j+S_j)=\Gamma(S_i,S_0+S_i)=1,\vspace{2mm}\\
x^{\lam_i}\cdot(x^{\lam_j}\cdot x^{\lam_k})=\Gamma(S_j,S_k)\Gamma(S_i,S_j+S_k)x^{\lam_i+\lam_j+\lam_k},\quad \vspace{2mm}\\
\;[L_{x^{\sg}},L_{x^{\tau}}](x^\gamma)=0\hbox{ if at least one of $\sg,\tau$ or $\gamma$ belongs to $S_0$},\vspace{2mm}\\
\hbox{$[L_{x^{\sg}},L_{x^{\tau}}]=0,\quad (\sg,\tau\in S_i)$},\vspace{2mm}\\
\hbox{$[L_{x^{\sg}}, L_{x^{\tau}}]=[L_{x^{\sg+\tau+\tau_j}}, L_{x^{-\tau_j}}]$},\quad(\sg,\tau\in S).
	\end{array}
	\end{equation}
Using these facts, we now claim that
\begin{equation}\label{eq:m1}
	[{x^{\lam_i+n\sg}}, 
	\bar{x}^{-\lam_i+n'\sg}]=\left\{\begin{array}{ll}
	L_{x^{t\sg}}&\hbox{if } i=0\hbox{ or }t\in 2\bbbz\\
	L_{x^{t\sg}}&\hbox{if } i\not=0,\;t\in 2\bbbz+1,\;\sg\in S_0\cup S_i\\
	\pm[L_{{x^{\lam_i+(t-1)\sg}}},
	L_{x^{-\lam_i+\sg}}]&\hbox{otherwise,}
	\end{array}\right.
\end{equation}
where $0\leq i\leq m$, $\lam_i\in S_i$, $\sg\in R^0\setminus\{0\}$ with $\lam_i+\sg\in S$, and
$n+n'=t$, $t,n,n'\in\bbbz$.
	
To see this, we note from (\ref{eq:khat4}) and (\ref{eq:rokn3}) that
\begin{eqnarray*}
0\not=[{x^{\lam_i+n\sg}}, 
	\bar{x}^{-\lam_i+n'\sg}]
	&=&
	{x^{\lam_i+n\sg}}\Delta\; 
	{\bar{x}^{-\lam_i+n'\sg}}\\
	&=&{L_{\underset{J_1}{\underbrace{x^{\lam_i+n\sg}\cdot x^{-\lam_i+n'\sg}}}}} +
	\stackrel{J_2}{\overbrace{[L_{{x^{\lam_i+n\sg}}},L_{x^{-\lam_i+n'\sg}}]}}.
	\end{eqnarray*}
Since $\sg\in R^0=S+S$, we have $\sg=\tau_j+\tau_k+2\lam$, $\lam\in\Lam$, for some $0\leq j,k\leq m$. Then as $\lam_i+\sg\in S$, it forces $\lam_i+\tau_j+\tau_k\in S$.
Note that
	$$J_1
	=\Gamma(\lam_i+n\sg,-\lam_i+n'\sg)x^{t\sg}=\Gamma(S_i+n(S_j+S_k), S_i+n'(S_j+S_k))x^{t\sg}.$$
		Now, if $i=0$, or $n+n'\in 2\bbbz$ or $\sg\in S_0\cup S_i$, then 
	at least one of $\lam_i+n\sg$, $-\lam_i+n'\sg$ belongs to $S_i$ or $S_0$, or both belong to the same coset, implying that $J_2=0$ by (\ref{eq:rokn3}), and $\Gamma(\lam_i+n\sg,-\lam_i+n'\sg)=1$. Thus
	$[{x^{\lam_i+n\sg}}, 
	\bar{x}^{-\lam_i+n'\sg}]
=L_{{x^{t\sg}}}.$
	 
	Assume next that $i\not=0$, $n\in 2\bbbz$, $n'\in2\bbbz+1$, and $\sg\not\in S_0\cup S_i$.
	The latter implies that $S_i+\sg\not=S_i$ and $\sg+S_i\not=S_0$ so
	$\Gamma(\lam_i+n\sg ,-\lam_i+n'\sg)=
	\Gamma(S_i, S_i+\sg)=0$, therefore $J_1=0$.  On the other hand
	\begin{eqnarray*}
	J_2&=&
	[L_{{x^{\lam_i+n\sg}}},L_{x^{-\lam_i+n'\sg}}]\\
	\hbox{(by (\ref{eq:rokn3})}	&=&
	[L_{{x^{\lam_i+(n+n'-1)\sg}}},L_{x^{-\lam_i+\sg}}]\\
&=&
	[L_{{x^{\lam_i+(t-1)\sg}}},L_{x^{-\lam_i+\sg}}].
	\end{eqnarray*}
	If $i\not=0$, $n$ is odd, $n'$ is even and $\sg\not\in S_0\cup S_i$, we get 
	$J_2=-
	[L_{{x^{\lam_i+(k-1)\sg}}},L_{x^{-\lam_i+\sg}}],$  by symmetry. This completes the argument
	 that (\ref{eq:m1}) holds.


 Now let $\a\in R^\times$, $\sg\in R^0$ and $\a+\sg\in R$. According to Corollary \ref{cor:12},  $\hat{E}_{\a,\sg}$ is an affine Lie subalgebra of $E$. We want to describe the root spaces of $\hat{E}_{\a,\sg}$, associated to non-zero isotropic roots. Each non-zero isotropic root is of the form $t\sg$ for some $t\not=0$.
 We have $\a=\dot\a+\lam_i$ for some $\lam_i\in S$. Since $\a+\sg\in R^\times$, we get
$\lam_i+\sg\in S$. Now for $n,n'\in\bbbz$, we have 
$E_{\a+n\sg}=\bbbk x^{\lam_i+n\sg}$, $ E_{-\a+n\sg}=\bbbk \bar{x}^{-\lam_i+n'\sg}$, and 
\begin{eqnarray*}
0\not=[E_{\a+n\sg}, E_{-\a+n'\sg}]=
	\bbbk[ {x^{\lam_i+n\sg}}, 
	{\bar{x}^{-\lam_i+n'\sg}}].
	\end{eqnarray*}
This together with (\ref{eq:m1}) shows that
$$(\hat{E}_{\a,\sg})_{t\sg}
=\left\{\begin{array}{ll}
	L_{x^{t\sg}}&\hbox{if } i=0\hbox{ or }t\in 2\bbbz\\
	L_{x^{t\sg}}&\hbox{if } i\not=0,\;t\in 2\bbbz+1,\;\sg\in S_0\cup S_i\\
	\pm[L_{{x^{\lam_i+(t-1)\sg}}},
	L_{x^{-\lam_i+\sg}}]&\hbox{otherwise.}
	\end{array}\right.
	$$
\end{exa}

\section{Filtration in extended affine Lie algebras}\setcounter{equation}{0}\label{sec:filter}
Let $(E,\fm,\hh)$ be an extended affine Lie algebra with root system $R$. In this section, we discuss existence of (ascending) filtrations for $R$ and $E$. We begin with a convention.
\begin{con}\label{con:1}
	By an inclusion $(E',\fm',\hh')\sub (E,\fm,\hh)$ of two extended affine Lie algebras, we mean $E'$ is a subalgebra of $E$,
	$\fm'={\fm}_{|_{E'}}$, and $\hh'\sub\hh$. Also we recall from Definition \ref{def:n1} that
we write $(R',\fm',\v')\sub (R,\fm,\v)$ to indicate that $R'$ is a subsystem of $R$.
\end{con}
 
Let $R_0\sub R_1\sub\cdots\sub R_n$ be a sequence of subsystems of $R$, where for each $i$, $R_i$ is closed in $R_{i+1}$. Let $\hh_0\sub\hh_1\sub\cdots
\sub\hh_n$ be some Lie covers for $R_0,\ldots, R_n$, respectively. Then by Proposition \ref{pro:khan2}, we get a filtered sequence of extended affine Lie algebras
$$E_{R_0,\hh_0}\sub\cdots\sub E_{R_n,\hh_n}.$$
If $R_n=R$ and $\hh_n=\hh$, we call $\{E_{R_i,\hh_i}\}_i$ a filtration for $E$.
 For different needs, appropriate filtrations may be applied.
We examine some here. 

\begin{parag}\label{parag:1}
	Let $(R,\fm,\v)$ be an extended affine root system. According to \pref{pt1}, we have
$$R=(S+S)\cup(\rds +S)\cup (\rdl+L),$$
where $\rd=\rds\cup\rdl\cup\{0\}$ is an irreducible reduced finite root system and $S$ and $L$ are two semilattices
satisfying (\ref{eq:inter1}). To refer to this description of $R$, we write $R=R(\rd,S,L)$.
We also recall from \cite[Chapter II.\S 4(b)]{AABGP97} that if $\rdl$ is of type $B_\ell$, then $\Lam=\la S\ra$ has a lattice decomposition $\Lam=\Lam_1\oplus\Lam_2$ with $S=S_1\oplus\Lam_2$ and $L=2\Lam_2\oplus S_2$, where $S_1$ and $S_2$ are semilattices in $\Lam_1$ and $\Lam_2$, respectively. 

We recall that if $S$ is a semilattice in $\Lam$, then $S=\cup_{i=0}^m(\tau_i+2\Lam)$ where $\tau_0=0$, $\tau_i$'s represent distinct costs of $2\Lam$ in $\Lam$, and $\La S\ra=\Lam$. In this case the integer $m$ is called the {\it index} of $S$, denoted $m=\ind(S)$. If $\ind(S)=\rank\,\Lam$, then $\tau_i$'s can be chosen such that $\{\tau_1,\ldots,\tau_m\}$ is a $\bbbz$-basis of $\Lam$.


To follow the next two results, we encourage the reader to consult \cite[Chapter II.\S1]{AABGP97} or \cite{Az97} for more details on semilattices. 
\end{parag}

\begin{lem}\label{lem:s1}
Let $S$ be the semilattice as in \pref{parag:1} and $\Lam=\la S\ra$. Let $\Lam=\Lam_1\oplus\Lam_2$, where $\Lam_1$ and $\Lam_2$ are two lattices with $\Lam_1\not=\{0\}$.
Let $S=S_1\oplus \Lam_2$ where $S_1$ is a semilattice in $\Lam_1$ with $\ind(S_1)=\rank\,\Lam_1$.
Assume that $\tilde{S}$ is a semilattice of the form $\tilde{S}=\tilde{S}_1\oplus\Lam'_2$, where $\tilde{S}_1\sub S_1$ is a semilattice with $(\tilde{S}_1\cap 2\Lam_1)+\tilde{S}_1\sub\tilde{S}_1$  and $\Lam'_2\sub\Lam_2$ is a lattice. 
Then $\la \tilde{S}\ra\cap R^0= \tilde{S}+\tilde{S}$. 
\end{lem} 

\proof We have $\tilde{S}+\tilde{S}\sub S+S=R^0$. Thus $\tilde{S}+\tilde{S}\sub\la \tilde{S}\ra\cap R^0$. 
We now show the reverse inclusion. Let $n=\rank\,\Lam_1$. Since
$\ind(S_1)=n$, we may assume $S_1=\cup_{i=0}^n(\sg_i+2\Lam_1)$ 
where $\sg_0=0$ and $\Lam_1=\sum_{i=1}^n\bbbz\sg_i$. Let $\eta\in\la \tilde{S}\ra\cap R^0$.
Since $\eta\in\la\tilde S\ra=\la\tilde{S}_1\ra \oplus\Lam'_2$, we have $\eta=\eta'+\lam_2'$ for some $\eta'\in\la \tilde{S}_1\ra$ and $\lam'_2\in\Lam'_2$. Since $\tilde{S}_1\sub S_1$, we have
$\tilde{S}_1=\cup_{i=0}^n(\tilde{S}_1\cap(\sg_i+2\Lam_1))$ and so 
$\eta'=\sum_{i=0}^n m_i(\sg_i+2\mu_i)$ for some $m_i\in\bbbz$ and $\mu_i\in\Lam_1$,  where $\sg_i+2\mu_i\in\tilde{S}_1$ if $m_i\not=0$. 
Since $\eta\in R^0=S+S=S_1+S_1+\Lam_2$, we have
$\eta=\sg_j+\sg_k++2\lam_1+\lam_2$ for some $0\leq j,k\leq n$, $\lam_1\in\Lam_1$ and $\lam_2\in\Lam_2$.
Then
\begin {eqnarray*}
\eta'-2m_0\mu_0&=&
\sum_{i=1}^n m_i(\sg_i+2\mu_i)\\
&=&
\eta-\lam'_2-2m_0\mu_0\\
&=&
\sg_j+\sg_k+2\lam_1+\lam_2-\lam'_2-2m_0\mu_0.
\end{eqnarray*} 
	But $\eta'-2m_0\mu_0\in\Lam_1$ and $\lam_2-\lam'_2\in\Lam_2+\Lam'_2\sub\Lam_2$. This forces
$$\sum_{i=1}^n m_i(\sg_i+2\mu_i)=\sg_j+\sg_k+2\lam_1,\hbox{ for some }\lam_1\in\Lam_1,\qquad (\star).$$
Now if $1\leq j\not= k\leq n$, then $(\star)$ implies $m_i=2m'_i$ for $i\not\in\{0,j,k\}$, and so we have
\begin{eqnarray*}
	\eta&=&m_j(\sg_j+2\mu_j)+m_k(\sg_k+2\mu_k)+2m_0\mu_0\\
	&&\qquad\qquad+
	\sum_{\{i\mid i\not=0,j,k\}}2m'_i(\sg_i+2\mu_i)+\lam'_2\\
&\in& 
\tilde{S_1}+\tilde{S}_1+(\tilde{S}_1\cap 2\Lam_1)+2\la \tilde{S}_1\ra+\Lam'_2\\
&\sub&\tilde{S}_1+\tilde{S}_1+2\la\tilde{S}_1\ra+\Lam'_2\\ 
&=&\tilde{S}_1+\tilde{S}_1+\Lam'_2=\tilde{S}+\tilde{S},
\end{eqnarray*}
where the inclusion $``\sub"$ follows from our assumption 
$\tilde{S}_1+(\tilde{S}_1\cap 2\Lam_1)\sub\tilde{S}_1$.

If $j=k$, then from $(\star)$ we have 
$\eta'=2m_0\mu_0+2\sg_j$ mod $2\Lam_1$. Thus $\eta'=2m_0\mu_0+\sum_{i=1}^n2m'_i(\sg_i+2\mu_i)$ for some $m'_i\in\bbbz$. This gives
$$\eta=\eta'+\lam'_2\in\tilde{S}_1+2\la\tilde{S}_1\ra+
\Lam'_2\sub\tilde{S_1}+
\Lam'_2\sub\tilde{S}+\tilde{S}.$$

The only remaining possibility is that only one of $j,k$ is $0$, say
$j=0$ and $1\leq k\leq n$. Then by $(\star)$ we have $m_i=2m'_i$ for $i\not\in\{0,k\}$, where $m'_i\in\bbbz$.
Thus
$$\eta=\eta'+\lam'_2=2m_0\mu_0+\sum_{\{i\mid i\not=0,k\}}2m'_i(\sg_i+2\mu_i)
\in\tilde{S}_1+2\la\tilde{S}_1\ra+\Lam'_2\sub\tilde{S}+\tilde{S}.$$
This completes the proof.\qed

If $R_1$ is a subsystem of $R$,  we say $R_1$ is a {\it canonical} subsystem if
$R_1=R(\rd,S_1,L_1)$. More precisely, $R_1$ is canonical if
$$R_1=(S_1+S_1)\cup (\rds+S_1)\cup (\rdl+L_1).$$
Therefore, a canonical subsystem has the same type and rank of $R$, but it may have smaller nullity. 
We now consider a general construction of canonical subsystems. Recall from \pref{parag:n1}, the subset $\tilde{R}_T$ of $R$ associated to a subset $T$ of $R^\times$.

\begin{lem}\label{lem:10}
	Let $U$ be a subgroup of $\la R^0\ra$, $T=(\dot R+U)\cap R^\times,$ and $\tilde S=S\cap U$. Then $T=(\dot R+\tilde S)\cap R^\times=
	\tilde{R}^\times_T$ and $\tilde{R}^0_T=\tilde S+\tilde S$.
	\end{lem}
	
	\proof
	From \pref{pt1}, we observe that $R^\times=(\dot R+S)\cap R^\times$. Thus $T=(\dot R+\tilde S)\cap R^\times$, where $\tilde S=S\cap U$. 
	
	Next, let $\a\in \tilde{R}_T^\times=\la T\ra\cap R^\times.$ Since $\a\in R^\times\sub\dot R+S$, we have $\a=\dot\a+\sg$ for some
	$\dot\a\in\dot{R}^\times$, $\sg\in S$. Since $\a\in\la T\ra$, we have $\a=\dot\a+\sg\in\la\dot R\ra\oplus A$, forcing that $\sg\in U$. Therefore,
	$\a\in\dot R+\tilde S$ and so
	$$\tilde{R}^\times_T\sub (\dot R+\tilde S)\cap R^\times=T\sub {\tilde R}^\times_T,$$
	proving the first assertion.

Finally,
\begin{eqnarray*}
	{\tilde R}_T^0&=&({\tilde R}^\times_T-{\tilde R}^\times_T)\cap\v^0\\
	&=&
\big(
((\dot R+\tilde S)\cap R^\times)-((\dot R+\tilde S)\cap R^\times)
\big)
\cap\v^0\\
	&=&\tilde S+\tilde S,
	\end{eqnarray*}
	and we are done.\qed

\begin{pro}\label{pro:3}
	Let $R=R(\rd,S,L)$ be an extended affine root system, 
$U$ be a subgroup of $\la R^0\ra$ and $T=(\rd+U)\cap R^\times$.  Then $\tilde{R}_T$ is a closed canonical subsystem of $R$ provided that one of the following conditions holds:
	 
	 (1) $\dot R\not= A_1,\;B_\ell$.

(2) $\dot R=A_1$, and either $S$ is a lattice or $ind(S)=\rank\, \Lam$.
	
	(3)  $\dot R=B_\ell$, and either $S$ is a lattice or $\ind(S_1)=\rank\,\Lam_1$.
	\end{pro}

	\proof
	Since $\rd^\times$ is connected, $T$ is also connected. By \pref{parag:n1}, the subset $\tilde{R}_T$ is a real-closed subsystem of $R$ and is clearly canonical. Therefore,
	to complete the proof, it remains to prove that if $\eta_1,\eta_2\in\tilde{R}_T^0$ and
	$\eta_1+\eta_2\in R^0$, then $\eta_1+\eta_2\in\tilde{R}_T^0$. By Lemma \ref{lem:10}, we have
	$$\tilde{R}_T=(\tilde S+\tilde S)\cup ((\dot R+\tilde S)\cap R^\times)$$
	 where $\tilde S=S\cap U$. In particular,  $\tilde{R}_T$ has the same type as $R$. 
	
	If $S$ is a lattice, then $\tilde{S}$ and hence $\tilde{R}^0_T$ are also lattices, and the claim follows immediately. Therefore, by \cite[Construction II.2.32, and Theorem II.2.37]{AABGP97}, the statement holds for all types other than $A_1$, and $B_\ell$, and also for these two types when $S$ is a lattice.

Next, assume that $R$ has type $A_1$ and  $\ind(S)=\rank\,\Lam$. Since $S+2\Lam=S$, we have
$$(\tilde{S}\cap 2\Lam)+\tilde S=((S\cap U)\cap 2\Lam)+(S\cap U)=(U\cap 2\Lam)+(S\cap U)\sub S\cap U=\tilde{S}.$$
Therefor, the conditions of Lemma \ref{lem:s1} are satisfied with $\Lam_1=\Lam=\la S\ra$ and $\Lam_2=\{0\}$. Thus
	$$\eta_1+\eta_2\in\la\tilde{R}_T^0\ra\cap R^0=\la\tilde{S}\ra\cap R^0=\tilde{S}+\tilde{S}=\tilde{R}_T^0.$$
	
	Finally, assume that $R$ is of type $B_\ell$ with  $\ind(S_1)=\rank\,\Lam_1$.
Set $\tilde{S}_1=S_1\cap A$ and $\Lam'_2=A\cap\Lam_2$. Then, using the assumptions in the statement, we obtain 
$$\tilde{S}=\tilde{S}_1\oplus\Lam'_2,\quad \tilde{S}_1\sub S_1,\quad \Lam'_2\sub\Lam_2,$$ and 
$$(\tilde{S}_1\cap 2\Lam_1)+\tilde{S}_1=
(A\cap S_1\cap 2\Lam_1)+(A\cap S_1)=A\cap 2\Lam_1+(A\cap S_1)\sub  A\cap S_1=\tilde{S}_1.$$
 Thus $S$, $S_1$, $\tilde{S}$, $\tilde{S}_1$ and $\Lam'_2$ satisfy the conditions of Lemma \ref{lem:s1}. Hence
	 $$\eta_1+\eta_2\in\la\tilde{R}^0_T\ra=\la\tilde S\ra\cap R^0=\tilde{S}+\tilde{S}=\tilde{R}_T^0.$$  
	\qed

\begin{pro}\label{pro:1} Let $(E,\fm,\hh)$ be an extended affine Lie algebra of nullity $\nu$, with root system $(R,\fm,\v)$.
Further, if $R$ is  of type $A_1$ or $B_\ell$ assume either that $S$ is a lattice, or otherwise that $\ind(S)=\rank\,\Lam$ in type $A_1$, and $\ind(S_1)=\rank\,\Lam_1$ in type $B_\ell$. 
Then $R$ admits a filtration of extended affine root systems
	$$(R_0,\fm_0,\v_0)\sub\cdots\sub(R_\nu,\fm_\nu,\v_\nu)=(R,\fm,\v),$$
	and a corresponding  filtration  of extended affine Lie algebras 
	$$(E_0,\fm_0,\hh_0)\sub\cdots\sub(E_\nu,\fm_\nu,\hh_\nu)=(E,\fm,\hh),$$
such that for $0\leq k\leq\nu$, the root system $R_k$ corresponds to the Lie algebra $E_k$, preserving the type and rank of $R$, and having nullity $k$.
	\end{pro}

	\proof
We have $\Lam=\sum_{i=1}^\nu\bbbz\sg_i$ where $\sg_i\in S$ for each $i$. For $0\leq k\leq\nu-1$ we set
	$$\sg_0=0,\quad U_k=\sum_{i=0}^k\bbbz\sg_i,\andd T_k=(\rd+U_k)\cap R^\times.$$
	By Proposition  \ref{pro:3}, for each $k$, the set $\tilde{R}_{T_k}$ is a canonically closed subsystem of $R$. 
	Now, by setting ${R}_i=\tilde{R}_{T_i}$, we have a filtration
	$$\dot R={R}_0\sub {R}_1\sub\cdots\sub {R}_{\nu-1}\sub {R}_\nu=R,$$
	for $R$, where each ${R}_k$ is an extended affine subsystem with the same type and rank as $R$, of nullity $k$.
	
	Next, as we discussed in Section \ref{sec:local}, using \cite[Theorem 11.11]{Rom05}, the orthogonal complement $\dot\hh^\perp$ of $\dot\hh$ in $\hh$ contains subspaces 
$$
		\begin{array}{c}
			\dim{\hat\hh}_k^0=\dim\sum_{i=0}^k\bbbk t_{\sg_i},\\
			\hbox{the form restricted to $\sum_{i=0}^k\bbbk t_{\sg_i}\oplus\hat\hh^0_k$ is non-degenerate},\\
			({\hat\hh}_k^0,{\hat\hh}_k^0)=\{0\}.
		\end{array}
		$$
	Moreover, the subspaces ${\hat\hh}^0_i$ can be chosen such that
	$$\{0\}={\hat\hh}^0_0\sub{\hat\hh}^0_1\sub\cdots\sub{\hat\hh}^0_{\nu-1}.$$
	Now, we set
	$$\hh_k=\dot\hh\oplus(\sum_{i=0}^k\bbbk t_{\sg_i})\oplus{\hat\hh}^0_k,\andd\hh_\nu=\hh,\qquad(0\leq k\leq\nu-1).$$

Next, using the notation of Section \ref{sec:local}, we set $(E_k,\fm_k,\hh_k):=E_{R_k,\hh_k}$.
By Proposition \ref{pro:khan2}, we get the required filtration as is claimed in the statement.\qed

\end{document}